
%
%
%
\def\unredoffs{} \def\redoffs{\voffset=-.31truein\hoffset=-.59truein}
\def\speclscape{\special{ps: landscape}}
%
%
%
%
\newbox\leftpage \newdimen\fullhsize \newdimen\hstitle \newdimen\hsbody
\tolerance=1000\hfuzz=2pt
\catcode`\@=11 
%
\ifx\answ\bigans\message{(This will come out unreduced.}
\magnification=1200\unredoffs\baselineskip=16pt plus 2pt minus 1pt
\hsbody=\hsize \hstitle=\hsize 
\else\message{(This will be reduced.} \let\l@r=L
\magnification=1000\baselineskip=16pt plus 2pt minus 1pt \vsize=7truein
\redoffs \hstitle=8truein\hsbody=4.75truein\fullhsize=10truein\hsize=\hsbody
\output={\ifnum\pageno=0 
  \shipout\vbox{\speclscape{\hsize\fullhsize\makeheadline}
    \hbox to \fullhsize{\hfill\pagebody\hfill}}\advancepageno
  \else
  \almostshipout{\leftline{\vbox{\pagebody\makefootline}}}\advancepageno
  \fi}
\def\almostshipout#1{\if L\l@r \count1=1 \message{[\the\count0.\the\count1]}
      \global\setbox\leftpage=#1 \global\let\l@r=R
 \else \count1=2
  \shipout\vbox{\speclscape{\hsize\fullhsize\makeheadline}
      \hbox to\fullhsize{\box\leftpage\hfil#1}}  \global\let\l@r=L\fi}
\fi
%
\newcount\yearltd\yearltd=\year\advance\yearltd by -1900

%

%
%

\def\draftmode{\message{ DRAFTMODE }\def\draftdate{{\rm preliminary draft:
\number\month/\number\day/\number\yearltd\ \ \hourmin}}%
\headline={\hfil\draftdate}\writelabels\baselineskip=20pt plus 2pt minus 2pt
 {\count255=\time\divide\count255 by 60 \xdef\hourmin{\number\count255}
  \multiply\count255 by-60\advance\count255 by\time
  \xdef\hourmin{\hourmin:\ifnum\count255<10 0\fi\the\count255}}}
\def\nolabels{\def\wrlabeL##1{}\def\eqlabeL##1{}\def\reflabeL##1{}}
\def\writelabels{\def\wrlabeL##1{\leavevmode\vadjust{\rlap{\smash%
{\line{{\escapechar=` \hfill\rlap{\sevenrm\hskip.03in\string##1}}}}}}}%
\def\eqlabeL##1{{\escapechar-1\rlap{\sevenrm\hskip.05in\string##1}}}%
\def\reflabeL##1{\noexpand\llap{\noexpand\sevenrm\string\string\string##1}}}
\nolabels
%
\global\newcount\secno \global\secno=0
\global\newcount\meqno \global\meqno=1
\def\newsec#1{\global\advance\secno by1\message{(\the\secno. #1)}
\global\subsecno=0\eqnres@t\noindent{\bf\the\secno. #1}
\writetoca{{\secsym} {#1}}\par\nobreak\medskip\nobreak}
\def\eqnres@t{\xdef\secsym{\the\secno.}\global\meqno=1\bigbreak\bigskip}
\def\sequentialequations{\def\eqnres@t{\bigbreak}}\xdef\secsym{}
\global\newcount\subsecno \global\subsecno=0
\def\subsec#1{\global\advance\subsecno by1\message{(\secsym\the\subsecno. #1)}
\ifnum\lastpenalty>9000\else\bigbreak\fi
\noindent{\it\secsym\the\subsecno. #1}\writetoca{\string\quad
{\secsym\the\subsecno.} {#1}}\par\nobreak\medskip\nobreak}
\def\appendix#1#2{\global\meqno=1\global\subsecno=0\xdef\secsym{\hbox{#1.}}
\bigbreak\bigskip\noindent{\bf Appendix #1. #2}\message{(#1. #2)}
\writetoca{Appendix {#1.} {#2}}\par\nobreak\medskip\nobreak}
%
%
\def\eqnn#1{\xdef #1{(\secsym\the\meqno)}\writedef{#1\leftbracket#1}%
\global\advance\meqno by1\wrlabeL#1}
\def\eqna#1{\xdef #1##1{\hbox{$(\secsym\the\meqno##1)$}}
\writedef{#1\numbersign1\leftbracket#1{\numbersign1}}%
\global\advance\meqno by1\wrlabeL{#1$\{\}$}}
\def\eqn#1#2{\xdef #1{(\secsym\the\meqno)}\writedef{#1\leftbracket#1}%
\global\advance\meqno by1$$#2\eqno#1\eqlabeL#1$$}
%
\newskip\footskip\footskip14pt plus 1pt minus 1pt 
\def\footnotefont{\ninepoint}\def\f@t#1{\footnotefont #1\@foot}
\def\f@@t{\baselineskip\footskip\bgroup\footnotefont\aftergroup\@foot\let\next}
\setbox\strutbox=\hbox{\vrule height9.5pt depth4.5pt width0pt}
\global\newcount\ftno \global\ftno=0
\def\foot{\global\advance\ftno by1\footnote{$^{\the\ftno}$}}
%
\newwrite\ftfile
\def\footend{\def\foot{\global\advance\ftno by1\chardef\wfile=\ftfile
$^{\the\ftno}$\ifnum\ftno=1\immediate\openout\ftfile=foots.tmp\fi%
\immediate\write\ftfile{\noexpand\smallskip%
\noexpand\item{f\the\ftno:\ }\pctsign}\findarg}%
\def\footatend{\vfill\eject\immediate\closeout\ftfile{\parindent=20pt
\centerline{\bf Footnotes}\nobreak\bigskip\input foots.tmp }}}
\def\footatend{}
%
%
\global\newcount\refno \global\refno=1
\newwrite\rfile
%
\def\ref{\nref}
\def\nref#1{\xdef#1{[\the\refno]}\writedef{#1\leftbracket#1}%
\ifnum\refno=1\immediate\openout\rfile=refs.tmp\fi
\global\advance\refno by1\chardef\wfile=\rfile\immediate
\write\rfile{\noexpand\item{#1\ }\reflabeL{#1\hskip.31in}\pctsign}\findarg}
\def\findarg#1#{\begingroup\obeylines\newlinechar=`\^^M\pass@rg}
{\obeylines\gdef\pass@rg#1{\writ@line\relax #1^^M\hbox{}^^M}%
\gdef\writ@line#1^^M{\expandafter\toks0\expandafter{\striprel@x #1}%
\edef\next{\the\toks0}\ifx\next\em@rk\let\next=\endgroup\else\ifx\next\empty%
\else\immediate\write\wfile{\the\toks0}\fi\let\next=\writ@line\fi\next\relax}}
\def\striprel@x#1{} \def\em@rk{\hbox{}}
\def\lref{\begingroup\obeylines\lr@f}
\def\lr@f#1#2{\gdef#1{\ref#1{#2}}\endgroup\unskip}

\def\addref#1{\immediate\write\rfile{\noexpand\item{}#1}} 
\def\startrefs#1{\immediate\openout\rfile=refs.tmp\refno=#1}
\def\refs#1{\count255=1[\r@fs #1{\hbox{}}]}
\def\r@fs#1{\ifx\und@fined#1\message{reflabel \string#1 is undefined.}%
\nref#1{need to supply reference \string#1.}\fi%
\vphantom{\hphantom{#1}}\edef\next{#1}\ifx\next\em@rk\def\next{}%
\else\ifx\next#1\ifodd\count255\relax\xref#1\count255=0\fi%
\else#1\count255=1\fi\let\next=\r@fs\fi\next}
%

%
\newwrite\ffile\global\newcount\figno \global\figno=1
\def\fig{fig.~\the\figno\nfig}
\def\nfig#1{\xdef#1{fig.~\the\figno}%
\writedef{#1\leftbracket fig.\noexpand~\the\figno}%
\ifnum\figno=1\immediate\openout\ffile=figs.tmp\fi\chardef\wfile=\ffile%
\immediate\write\ffile{\noexpand\medskip\noexpand\item{Fig.\ \the\figno. }
\reflabeL{#1\hskip.55in}\pctsign}\global\advance\figno by1\findarg}
\def\xfig{\expandafter\xf@g}\def\xf@g fig.\penalty\@M\ {}
\def\figs#1{figs.~\f@gs #1{\hbox{}}}
\def\f@gs#1{\edef\next{#1}\ifx\next\em@rk\def\next{}\else
\ifx\next#1\xfig #1\else#1\fi\let\next=\f@gs\fi\next}
\newwrite\lfile
{\escapechar-1\xdef\pctsign{\string\%}\xdef\leftbracket{\string\{}
\xdef\rightbracket{\string\}}\xdef\numbersign{\string\#}}

\def\writestop{\def\writestoppt{\immediate\write\lfile{\string\pageno%
\the\pageno\string\startrefs\leftbracket\the\refno\rightbracket%
\string\def\string\secsym\leftbracket\secsym\rightbracket%
\string\secno\the\secno\string\meqno\the\meqno}\immediate\closeout\lfile}}
\def\writestoppt{}\def\writedef#1{}
\def\seclab#1{\xdef #1{\the\secno}\writedef{#1\leftbracket#1}\wrlabeL{#1=#1}}
\def\subseclab#1{\xdef #1{\secsym\the\subsecno}%
\writedef{#1\leftbracket#1}\wrlabeL{#1=#1}}
\newwrite\tfile \def\writetoca#1{}
\def\leaderfill{\leaders\hbox to 1em{\hss.\hss}\hfill}
\def\writetoc{\immediate\openout\tfile=toc.tmp
   \def\writetoca##1{{\edef\next{\write\tfile{\noindent ##1
   \string\leaderfill {\noexpand\number\pageno} \par}}\next}}}
%
%
%
\def\listtoc{\centerline{\bf Contents}\nobreak\medskip{\baselineskip=12pt
 \parskip=0pt\catcode`\@=11 \input toc.tmp \catcode`\@=12 \bigbreak\bigskip}}
\catcode`\@=12 
%
\edef\tfontsize{\ifx\answ\bigans scaled\magstep3\else scaled\magstep4\fi}
\font\titlerm=cmr10 \tfontsize \font\titlerms=cmr7 \tfontsize
\font\titlermss=cmr5 \tfontsize \font\titlei=cmmi10 \tfontsize
\font\titleis=cmmi7 \tfontsize \font\titleiss=cmmi5 \tfontsize
\font\titlesy=cmsy10 \tfontsize \font\titlesys=cmsy7 \tfontsize
\font\titlesyss=cmsy5 \tfontsize \font\titleit=cmti10 \tfontsize
\skewchar\titlei='177 \skewchar\titleis='177 \skewchar\titleiss='177
\skewchar\titlesy='60 \skewchar\titlesys='60 \skewchar\titlesyss='60
\def\titlefont{\def\rm{\fam0\titlerm}
\textfont0=\titlerm \scriptfont0=\titlerms \scriptscriptfont0=\titlermss
\textfont1=\titlei \scriptfont1=\titleis \scriptscriptfont1=\titleiss
\textfont2=\titlesy \scriptfont2=\titlesys \scriptscriptfont2=\titlesyss
\textfont\itfam=\titleit \def\it{\fam\itfam\titleit}\rm}
 \ifx\answ\bigans\else scaled\magstep1\fi
\ifx\answ\bigans\else

 \font\absi=cmmi10 scaled\magstep1
\font\absis=cmmi7 scaled\magstep1 \font\absiss=cmmi5 scaled\magstep1
\font\abssy=cmsy10 scaled\magstep1 \font\abssys=cmsy7 scaled\magstep1
\font\abssyss=cmsy5 scaled\magstep1 
\skewchar\absi='177 \skewchar\absis='177 \skewchar\absiss='177
\skewchar\abssy='60 \skewchar\abssys='60 \skewchar\abssyss='60
\fi
\font\ninerm=cmr9 \font\sixrm=cmr6 \font\ninei=cmmi9 \font\sixi=cmmi6
\font\ninesy=cmsy9 \font\sixsy=cmsy6 \font\ninebf=cmbx9
\font\nineit=cmti9 \font\ninesl=cmsl9 \skewchar\ninei='177
\skewchar\sixi='177 \skewchar\ninesy='60 \skewchar\sixsy='60
\def\ninepoint{\def\rm{\fam0\ninerm}
\textfont0=\ninerm \scriptfont0=\sixrm \scriptscriptfont0=\fiverm
\textfont1=\ninei \scriptfont1=\sixi \scriptscriptfont1=\fivei
\textfont2=\ninesy \scriptfont2=\sixsy \scriptscriptfont2=\fivesy
\textfont\itfam=\ninei \def\it{\fam\itfam\nineit}\def\sl{\fam\slfam\ninesl}%
\textfont\bffam=\ninebf \def\bf{\fam\bffam\ninebf}\rm}
%
%

\hyphenation{anom-aly anom-alies coun-ter-term coun-ter-terms}
\def\inv{^{\raise.15ex\hbox{${\scriptscriptstyle -}$}\kern-.05em 1}}

\def\Dsl{\,\raise.15ex\hbox{/}\mkern-13.5mu D} 
\def\dsl{\raise.15ex\hbox{/}\kern-.57em\partial}

\def\lspace{\ifx\answ\bigans{}\else\qquad\fi}
\def\lbspace{\ifx\answ\bigans{}\else\hskip-.2in\fi} 
\def\boxeqn#1{\vcenter{\vbox{\hrule\hbox{\vrule\kern3pt\vbox{\kern3pt
    \hbox{${\displaystyle #1}$}\kern3pt}\kern3pt\vrule}\hrule}}}
\def\mbox#1#2{\vcenter{\hrule \hbox{\vrule height#2in
        \kern#1in \vrule} \hrule}}  
%

\def\darr#1{\raise1.5ex\hbox{$\leftrightarrow$}\mkern-16.5mu #1}

\def\roughly#1{\raise.3ex\hbox{$#1$\kern-.75em\lower1ex\hbox{$\sim$}}}

%
%


\def\frac#1#2{{#1\over#2}}

\def\journal#1&#2(#3){\unskip, #1~\bf #2 \rm(19#3) }
\def\andjournal#1&#2(#3){\sl #1~\bf #2 \rm (19#3) }

\def\bra#1{\left\langle #1\right|}
\def\ket#1{\left| #1\right\rangle}

\catcode`\@=11\def\slash#1{\mathord{\mathpalette\c@ncel{#1}}}
\overfullrule=0pt
\def\steepslash{\c@ncel}
\def\frac#1#2{{#1\over #2}}

\def\:{\!:\!}
\def\inbar{\,\vrule height1.5ex width.4pt depth0pt}
\def\IQ{\relax\,\hbox{$\inbar\kern-.3em{\rm Q}$}}
\def\IB{\relax{\rm I\kern-.18em B}}
\def\IC{\relax\hbox{$\inbar\kern-.3em{\rm C}$}}
\def\IP{\relax{\rm I\kern-.18em P}}
\def\IR{\relax{\rm I\kern-.18em R}}
\def\ZZ{\relax\ifmmode\mathchoice
{\hbox{Z\kern-.4em Z}}{\hbox{Z\kern-.4em Z}}
{\lower.9pt\hbox{Z\kern-.4em Z}}
{\lower1.2pt\hbox{Z\kern-.4em Z}}\else{Z\kern-.4em Z}\fi}

\catcode`\@=12

\def\npb#1(#2)#3{{ Nucl. Phys. }{B#1} (#2) #3}
\def\plb#1(#2)#3{{ Phys. Lett. }{#1B} (#2) #3}
\def\pla#1(#2)#3{{ Phys. Lett. }{#1A} (#2) #3}
\def\prl#1(#2)#3{{ Phys. Rev. Lett. }{#1} (#2) #3}
\def\mpla#1(#2)#3{{ Mod. Phys. Lett. }{A#1} (#2) #3}
\def\ijmpa#1(#2)#3{{ Int. J. Mod. Phys. }{A#1} (#2) #3}
\def\cmp#1(#2)#3{{ Comm. Math. Phys. }{#1} (#2) #3}
\def\cqg#1(#2)#3{{ Class. Quantum Grav. }{#1} (#2) #3}
\def\jmp#1(#2)#3{{ J. Math. Phys. }{#1} (#2) #3}
\def\anp#1(#2)#3{{ Ann. Phys. }{#1} (#2) #3}
\def\prd#1(#2)#3{{ Phys. Rev. } {D{#1}} (#2) #3}
\def\ptp#1(#2)#3{{ Progr. Theor. Phys. }{#1} (#2) #3}
\def\aom#1(#2)#3{{ Ann. Math. }{#1} (#2) #3}

\def\bs{\bigskip}

\def\bra{\langle}
\def\ket{\rangle}

\def\C{{\bf C}}

\def\F{{\bf F}}

\def\H{{\bf H}}

\def\L{{\bf L}}

\def\P{{\bf P}}

\def\R{{\bf R}}

\def\Z{{\bf Z}}
\def\cA{{\cal A}}
\def\cB{{\cal B}}
\def\cC{{\cal C}}

\def\cE{{\cal E}}
\def\cF{{\cal F}}

\def\cI{{\cal I}}

\def\cL{{\cal L}}

\def\cQ{{\cal Q}}

\def\cS{{\cal S}}

\def\cW{{\cal W}}

\input amssym

\def\gg{{\goth g}}
\def\gh{{\goth h}}

\def\gk{{\goth k}}

\def\gn{{\goth n}}

\def\gs{{\goth s}}
\def\gt{{\goth t}}

\def\cicy#1(#2|#3)#4{\left(\matrix{#2}\right|\!\!
                     \left|\matrix{#3}\right)^{{#4}}_{#1}}

\def\emptyset{\varnothing}

\def\ra{\rightarrow}

\def\bs{\bigskip}

\def\Box{{\,\lower0.9pt\vbox{\hrule
\hbox{\vrule height 0.2 cm \hskip 0.2 cm
\vrule height 0.2 cm}\hrule}\,}}

\global\newcount\thmno \global\thmno=0
\def\definition#1{\global\advance\thmno by1
\bigskip\noindent{\bf Definition \secsym\the\thmno. }{\it #1}
\par\nobreak\medskip\nobreak}
\def\question#1{\global\advance\thmno by1
\bigskip\noindent{\bf Question \secsym\the\thmno. }{\it #1}
\par\nobreak\medskip\nobreak}
\def\theorem#1{\global\advance\thmno by1
\bigskip\noindent{\bf Theorem \secsym\the\thmno. }{\it #1}
\par\nobreak\medskip\nobreak}
\def\proposition#1{\global\advance\thmno by1
\bigskip\noindent{\bf Proposition \secsym\the\thmno. }{\it #1}
\par\nobreak\medskip\nobreak}
\def\corollary#1{\global\advance\thmno by1
\bigskip\noindent{\bf Corollary \secsym\the\thmno. }{\it #1}
\par\nobreak\medskip\nobreak}
\def\lemma#1{\global\advance\thmno by1
\bigskip\noindent{\bf Lemma \secsym\the\thmno. }{\it #1}
\par\nobreak\medskip\nobreak}
\def\conjecture#1{\global\advance\thmno by1
\bigskip\noindent{\bf Conjecture \secsym\the\thmno. }{\it #1}
\par\nobreak\medskip\nobreak}
\def\exercise#1{\global\advance\thmno by1
\bigskip\noindent{\bf Exercise \secsym\the\thmno. }{\it #1}
\par\nobreak\medskip\nobreak}
\def\remark#1{\global\advance\thmno by1
\bigskip\noindent{\bf Remark \secsym\the\thmno. }{\it #1}
\par\nobreak\medskip\nobreak}
\def\problem#1{\global\advance\thmno by1
\bigskip\noindent{\bf Problem \secsym\the\thmno. }{\it #1}
\par\nobreak\medskip\nobreak}
\def\others#1#2{\global\advance\thmno by1
\bigskip\noindent{\bf #1 \secsym\the\thmno. }{\it #2}
\par\nobreak\medskip\nobreak}
\def\proof{\noindent Proof: }

\def\thmlab#1{\xdef #1{\secsym\the\thmno}\writedef{#1\leftbracket#1}\wrlabeL{#1=#1}}
%
%
\def\newsec#1{\global\advance\secno by1\message{(\the\secno. #1)}
\global\subsecno=0\thmno=0\eqnres@t\noindent{\bf\the\secno. #1}
\writetoca{{\secsym} {#1}}\par\nobreak\medskip\nobreak}
\def\eqnres@t{\xdef\secsym{\the\secno.}\global\meqno=1\bigbreak\bigskip}
\def\sequentialequations{\def\eqnres@t{\bigbreak}}\xdef\secsym{}
%

%
\newcount{\exnum}
\def\prob{\advance\exnum by 1
\bigskip\item{\the\exnum.}\ }
\newcount{\exnum}
\def\next{\advance\exnum by 1
\bigskip\noindent{\the\exnum.}\ }
\def\np{\vfill\eject}

\nopagenumbers

\ref\BHS{D. Ben-Zvi, R. Heluani, and M. Szczesny, Supersymmetry of the Chiral de Rham complex, Compos. Math. 144 (2) (2008) 503-521.}
\ref\CartanI{ H. Cartan, Notions d'alg\`ebre diff\'erentielle; application aux groupes de Lie et aux vari\'eti\'es o\`u op\`ere un groupe de Lie, Colloque de Topologie, C.B.R.M., Bruxelles 15-27 (1950).}
\ref\CartanII{ H. Cartan, La Transgression dans un groupe de Lie
et dans un espace fibr\'e principal, Colloque de Topologie,
C.B.R.M., Bruxelles 57-71 (1950).} 
\ref\DKV{M. Duflo, S. Kumar, and M. Vergne, Sur la Cohomologie \'Equivariante des Vari\'et\'es Diff\'erentiables, Ast\'erisque 215 (1993).}
\ref\F{ B. Feigin,  The Semi-Infinite Homology of Lie, Kac-Moody and Virasoro Algebras, Russian Math. Survey 39 (1984) no. 2 195-196.}
\ref\FF{ B. Feigin and E. Frenkel,  Semi-infinite Weil complex and the Virasoro algebra,  Comm. Math. Phys. 137 (1991), 617-639.}
\ref\FGZ{ I.B. Frenkel, H. Garland, and G.J. Zuckerman,  Semi-Infinite Cohomology and String Theory, Proc. Natl. Acad. Sci. USA Vol. 83, No. 22 (1986) 8442-8446.} 
\ref\FS{E. Frenkel and M. Szczesny, Chiral de Rham complex and orbifolds,  J. Algebraic Geom. 16 (2007), no. 4, 599--624.}
\ref\GKM{M. Goresky, R. Kottwitz, and R. MacPherson, Equivariant cohomology, Koszul duality, and the localization theorem, Invent. Math. 131, 25-83 (1998).}
 \ref\GMS{V. Gorbounov, F. Malikov, and V. Schechtman, Gerbes of chiral differential operators, Math. Res. Lett. 7 (2000), no. 1, 55--66.} 
 \ref\GS{ V. Guillemin and S. Sternberg, Supersymmetry and
Equivariant de Rham Theory, Springer, 1999.}
\ref\LLI{B. Lian and A. Linshaw, Chiral equivariant cohomology I, Adv. Math. 209, 99-161 (2007).}
\ref\LLSI{B. Lian, A. Linshaw, and B. Song, Chiral equivariant cohomology II, Trans. Am. Math. Soc. 360 (2008), 4739-4776.}
 \ref\MS{F. Malikov and V. Schechtman, Chiral de Rham complex, II, Differential topology, infinite-dimensional Lie algebras, and applications, 149--188, Amer. Math. Soc. Transl. Ser. 2, 194, Amer. Math. Soc., Providence, RI, 1999.}
\ref\MSV{F. Malikov, V. Schechtman, and A. Vaintrob, Chiral de Rham complex, Comm. Math. Phys. 204 (1999), no. 2, 439--473.}
\ref\OI{R. Oliver, Compact Lie group actions on disks, Math. Z. 149, 79-97 (1976).}
\ref\OII{R. Oliver, Fixed points of disk actions, Bull. Am. Mat. Soc, Vol.85, No. 2, 279-280 (1976).}
\ref\TD{T. tom Dieck, Transformation Groups, de Gruyter Studies in Mathematics 8 (1987).}
\ref\Witten{E. Witten, Two-dimensional models with (0,2) supersymmetry: perturbative aspects, Adv. Theor. Math. Phys. 11 (2007), no. 1, 1--63.}
\ref\Vor{A. Voronov, Semi-infinite homological algebra. Invent. Math. 113, 103-146 (1993).}

\centerline{\titlefont Chiral Equivariant Cohomology III}

\bs
\centerline{Bong H. Lian, Andrew R. Linshaw and Bailin Song}
\bs

\centerline{ \it{Dedicated to the memory of our friend and colleague Jerome P. Levine}}
\bs

\baselineskip=13pt plus 2pt minus 2pt
ABSTRACT.  This is the third of a series of papers on a new equivariant cohomology that takes values 
in a vertex algebra, and contains and generalizes the classical equivariant cohomology of a manifold with a Lie group action \` a la H. Cartan. In this paper, we compute this cohomology for spheres and show that for any simple connected group $G$, there is a sphere with infinitely many actions of $G$ which have distinct chiral equivariant cohomology, but identical classical equivariant cohomology. Unlike the classical case, the description of the chiral equivariant cohomology of spheres requires a substantial amount of new structural theory, which we fully develop in this paper. This includes a quasi-conformal structure, equivariant homotopy invariance, and the values of this cohomology on homogeneous spaces. These results rely on crucial features of the underlying vertex algebra valued complex that have no classical analogues.

\baselineskip=15pt plus 2pt minus 1pt
\parskip=\baselineskip


\def\listtoc{\centerline{\bf Contents}\nobreak\medskip{\baselineskip=12pt
 \parskip=0pt\catcode`\@=11
\noindent {1.} {Introduction} \leaderfill{2} \par 
\noindent \quad{1.1.} {Outline of main results} \leaderfill{6} \par 
\noindent \quad{1.2.} {General remarks about group actions on manifolds} \leaderfill{9} \par 
\noindent {2.} {Mayer-Vietoris Sequences} \leaderfill{10} \par 
\noindent {3.} {Homotopy Invariance of ${\bf H}^*_G({\cal Q}'(M))$} \leaderfill{13} \par 
\noindent \quad{3.1.} {Chiral chain homotopies} \leaderfill{13} \par 
\noindent \quad{3.2.} {Relative derivations} \leaderfill{16} \par 
\noindent {4.} {A Quasi-conformal Structure on ${\bf H}^*_G({\cal Q}'(M))$ and ${\bf H}^*_G({\cal Q}(M))$} \leaderfill{18} \par 
\noindent \quad{4.1.} {A vanishing criterion} \leaderfill{20} \par 
\noindent {5.} {${\bf H}^*_G({\cal Q}'(G/H))$ for Homogeneous Spaces $G/H$} \leaderfill{21} \par 
\noindent \quad{5.1.} {Finite-dimensionality of ${\bf H}^*_G({\cal Q}'(M))$ for compact $M$} \leaderfill{28} \par 
\noindent {6.} {The Structure of ${\bf H}^*_G({\cal Q}(M))$} \leaderfill{30} \par 
\noindent \quad{6.1.} {The ideal property of ${\bf H}^*_G({\cal Q}(M))_+$ and ${\bf H}^*_G({\cal Q}'(M))_+$} \leaderfill{31} \par 
\noindent {7.} {The Structure of ${\bf H}^*_G({\cal Q}'(M))$} \leaderfill{32} \par 
\noindent \quad{7.1.} {The case where $G$ is simple} \leaderfill{33} \par 
\noindent \quad{7.2.} {Chiral equivariant cohomology of spheres} \leaderfill{34} \par 
\noindent \quad{7.3.} {The case $G=G_1\times G_2$, where $G_1,G_2$ are simple} \leaderfill{38} \par 
\noindent \quad{7.4.} {The case where $G$ is a torus $T$} \leaderfill{42} \par 
\noindent \quad{7.5.} {The case $T=S^1\times S^1$ and $M = {\bf C}{\bf P}^2$} \leaderfill{44} \par 
\noindent {8.} {Concluding Remarks and Open Questions} \leaderfill{48} \par 

  \catcode`\@=12 \bigbreak\bigskip}}
\catcode`\@=12 

\listtoc

\np
\headline{\ifodd\pageno\rightheadline\else\leftheadline\fi}
\def\rightheadline{\tenrm\hfil Chiral Equivariant Cohomology III
\hfil\folio}
\def\leftheadline{\tenrm\folio\hfil B.H. Lian, A.R. Linshaw \& B. Song\hfil}

\newsec{Introduction}

Let $G$ be a compact Lie group with complexified Lie algebra $\gg$. For a topological $G$-space $M$, the equivariant cohomology $H^*_G(M)$ is defined to be $H^*((M\times \cE)/G))$, where $\cE$ is any contractible space on which $G$ acts freely. When $M$ is a manifold on which $G$ acts by diffeomorphisms, there is a de Rham model of $H^*_G(M)$ due to H. Cartan \CartanI\CartanII, and developed further by Duflo-Kumar-Vergne \DKV~and Guillemin-Sternberg \GS. In fact, one can define the equivariant cohomology $H^*_G(A)$ of any $G^*$-algebra
$A$. Taking $A$ to be the algebra $\Omega(M)$ of differential forms on $M$ gives us the de Rham model of $H^*_G(M)$, and $H^*_G(\Omega(M)) \cong H^*_G(M)$ by an equivariant version of the de Rham theorem.

In \LLI, the chiral equivariant cohomology $\H^*_G(\cA)$ of an $O(\gs\gg)$-algebra $\cA$ was introduced as a vertex algebra analogue of the equivariant cohomology of $G^*$-algebras. Examples of $O(\gs\gg)$-algebras include the semi-infinite Weil complex $\cW(\gg)$, which was introduced in \FF, and the chiral de Rham complex $\cQ(M)$ of a smooth $G$-manifold $M$. The chiral de Rham complex was introduced in \MSV, and has been studied from several different points of view in recent years \MS\GMS\FS\BHS\Witten. In \LLSI, the chiral equivariant cohomology functor was extended to the larger categories of $\gs\gg[t]$-algebras and $\gs\gg[t]$-modules. Our main example of an $\gs\gg[t]$-algebra which is {\it not} an $O(\gs\gg)$-algebra is the subalgebra $\cQ'(M)\subset \cQ(M)$ generated by the weight zero subspace. 
$\H^*_G(\cQ(M))$ and $\H^*_G(\cQ'(M))$ are both \lq\lq chiralizations" of $H^*_G(M)$, that is, vertex algebras equipped with weight gradings
$$\H^*_G(\cQ(M)) = \bigoplus_{m\geq 0} \H^*_G(\cQ(M))[m],~~~~\H^*_G(\cQ'(M)) = \bigoplus_{m\geq 0}\H^*_G(\cQ'(M))[m],$$ such that $\H^*_G(\cQ(M))[0] = H^*_G(M) = \H^*_G(\cQ'(M))[0]$. In the case $M = pt$, $\cQ(M) = \cQ'(M) = \C$, and $\H^*_G(\C)$ plays the role of $H^*_G(pt) = S(\gg^*)^G$ in the classical theory.

We briefly recall these constructions, following the notation in \LLI \LLSI. We will assume that the reader is familiar with the basic notions in vertex algebra theory. For a list of references, see page 102 of \LLI. A {\it differential vertex algebra} (DVA) is a degree graded vertex algebra $\cA^*=\oplus_{p\in\Z}\cA^p$ equipped with a vertex algebra derivation $d$ of degree 1 such that $d^2=0$. A DVA will be called {\it degree-weight graded} if it has an additional $\Z_{\geq 0}$-grading by weight, which is compatible with the degree in the sense that $\cA^p=\oplus_{n\geq0}\cA^p[n]$. There is an auxiliary structure on a DVA which is analogous to the structure of a $G^*$-algebra in \GS. 
Associated to $\gg$ is a Lie superalgebra $\gs\gg := \gg\triangleright \gg^{-1}$ with bracket $[(\xi,\eta),(x,y)]=([\xi,x],[\xi,y]-[x,\eta])$, which is equipped with a differential $d:(\xi,\eta)\mapsto(\eta,0)$. 
This differential extends to the loop algebra $\gs\gg[t,t^{-1}]$, and gives rise to a vertex algebra derivation on the corresponding current algebra $O(\gs\gg):= O(\gs\gg,0)$. Here $0$ denotes the zero bilinear form on $\gs\gg$. 

An $O(\gs\gg)$-algebra is a degree-weight graded DVA $\cA$ equipped with a DVA homomorphism $\rho:O(\gs\gg)\ra \cA$, which we denote by $(\xi,\eta)\ra L_{\xi} + \iota_{\eta}$. Although this definition makes sense for any Lie algebra $\gg$, we will assume throughout this paper that $\gg$ is the complexified Lie algebra of $G$, and we require $\cA$ to admit an action $\hat{\rho}:G\ra Aut(\cA)$ of $G$ by vertex algebra automorphisms which is compatible with the $O(\gs\gg)$-structure in the following sense:
\eqn\gsI{\frac{d}{dt}\hat{\rho}(exp(t\xi))|_{t=0} = L_{\xi}(0),}
\eqn\gsII{\hat{\rho} (g) L_{\xi}(n) \hat{\rho}(g^{-1}) = L_{Ad(g)(\xi)}(n),}
\eqn\gsIII{\hat{\rho}(g) \iota_{\xi}(n)\hat{\rho}(g^{-1}) =\iota_{Ad(g)(\xi)}(n),}
\eqn\gsIV{\hat{\rho}(g) d \hat{\rho} (g^{-1}) = d,}
for all $\xi\in\gg$, $g\in G$, and $n\in\Z$. These conditions are analogous to Equations (2.23)-(2.26) of \GS. In order for \gsI~to make sense, we must be able to differentiate along appropriate curves in $\cA$, which is the case in our main example $\cA=\cQ(M)$. 

In \LLSI, we observed that the chiral equivariant cohomology functor can be defined on the larger class of spaces which carry only a representation of the Lie subalgebra $\gs\gg[t]$ of $\gs\gg[t,t^{-1}]$.

\definition{An $\gs\gg[t]$-module is a degree-weight graded complex $(\cA,d_\cA)$ equipped with a Lie algebra homomorphism $\rho:\gs\gg[t]\ra End~\cA$, such that  for all $x\in\gs\gg[t]$ we have
\item{$\bullet$} $\rho(dx)=[d_\cA,\rho(x)]$
\item{$\bullet$} $\rho(x)$ has degree 0 whenever $x$ is even in $\gs\gg[t]$, and degree -1 whenever $x$ is odd, and has weight $-n$ if $x\in\gs\gg t^n$.}

As above, we will always assume that $\gg$ is the complexified Lie algebra of $G$, and $\cA$ has an action of $G$ which is compatible with the $\gs\gg[t]$-structure, i.e., \gsI-\gsIV~hold for $n\geq 0$. 

\definition{Given an $\gs\gg[t]$-module $(\cA,d)$, we define the chiral horizontal, invariant and basic subspaces of $\cA$ to be respectively
$$\eqalign{
\cA_{hor}&=\{a\in\cA|\rho(x)a=0~\forall x\in\gg^{-1}[t]\}\cr
\cA_{inv}&=\{a\in\cA|\rho(x)a=0~\forall x\in\gg[t],~\hat{\rho}(g)(a) = a~\forall g\in G\}\cr
\cA_{bas}&=\cA_{hor}\cap\cA_{inv}.
}$$
An $\gs\gg[t]$-module $(\cA,d)$ is called an $\gs\gg[t]$-algebra if it is also a DVA such that $\cA_{hor},\cA_{inv}$ are both vertex subalgebras of $\cA$, and $G$ acts by DVA automorphisms.}

When we are working with multiple groups $G,H ,\dots$ we will use the notations $\cA_{G-hor}$, $\cA_{G-inv}$, $\cA_{G-bas}$, etc., to avoid confusion. Given an $\gs\gg[t]$-module $(\cA,d)$, $\cA_{inv}$, $\cA_{bas}$ are both subcomplexes of $\cA$, but $\cA_{hor}$ is not a subcomplex of $\cA$ in general.
The Lie algebra $\gs\gg[t]$ is not required to act by derivations on a DVA $\cA$ to make it an $\gs\gg[t]$-algebra. If $(\cA,d)$ is an $O(\gs\gg)$-algebra, any subDVA $\cB$ which is closed under the operators $(L_\xi+\iota_\eta)\circ_p$, $p\geq0$, is an $\gs\gg[t]$-algebra.

\definition{For any $\gs\gg[t]$-module $(\cA,d_{\cA})$, we define its chiral basic cohomology $\H^*_{bas}(\cA)$ to be $ H^*(\cA_{bas},d_{\cA})$. We define its chiral equivariant cohomology $\H^*_{G}(\cA)$ to be $\H^*_{bas}(\cW(\gg)\otimes\cA)$. The differential on $\cW(\gg)\otimes \cA$ is $d_{\cW} \otimes 1 + 1\otimes d_{\cA}$, where $d_{\cW}=J(0)+K(0)$, as in \LLI. }

In this paper, our main focus is on the cases $\cA=\cQ(M)$ and $\cA = \cQ'(M)$, for a smooth $G$-manifold $M$. Recall from \LLSI~that for each $m\geq 0$, $\cQ_M[m]$ is a sheaf of vector spaces on $M$, and $\cQ(M)$ is the space of global sections of the weak sheaf of vertex algebras $\cQ_M=\oplus_{m\geq 0} \cQ_M[m]$ on $M$. Similarly, $\cQ'_M[m]$ is the subsheaf of $\cQ_M[m]$ generated by the weight zero subspace, and $\cQ'(M)$ is the space of global sections of the weak sheaf of abelian vertex algebras $\cQ'_M = \oplus_{m\geq 0} \cQ'_M[m]$. In this terminology, a {\it weak sheaf} is a presheaf which satisfies a slightly weaker version of the reconstruction axiom: 
\eqn\recon{0\ra\cF(U)\rightarrow\prod_i\cF(U_i)\rightrightarrows\prod_{i,j}\cF(U_i\cap U_j)} is exact for {\it finite} open covers $\{U_i\}$ of an open set $U$ (see Section 1.1 of \LLSI). Whenever we need to construct a global section of $\cQ_M$ or $\cQ'_M$ by gluing together local sections, these sections are always homogeneous of finite weight, so we may work inside the sheaf $\cQ_M[m]$ or $\cQ'_M[m]$ for some $m$.

In \LLI \LLSI~we proved a number of structural results about the chiral equivariant cohomology, which we recall below. 

\item{$\bullet$} The functor $\H^*_G(\cQ'(-))$ is contravariant in $M$ for any $G$. For fixed $M$, $\H^*_{(-)}(\cQ'(M))$ is not functorial in $G$ in general, but is contravariant with respect to abelian groups.

\item{$\bullet$} For any group $G$ of positive dimension, $\H^*_G(\C)$ contains nonzero classes in every positive weight. If $G$ is semisimple, $\H^*_G(\C)$ is a conformal vertex algebra with Virasoro element $\L$ of central charge zero.

\item{$\bullet$} For any $\gs\gg[t]$-algebra $\cA$, the canonical map
$$ \kappa_G:\H^*_{bas}(\cW(\gg))=\H^*_G(\C)\ra\H^*_G(\cA)$$ 
induced by $\cW(\gg)\hookrightarrow\cW(\gg)\otimes\cA$ is called the {\it chiral Chern-Weil map} of $\cA$. For $\cA = \cQ(M)$ or $\cA = \cQ'(M)$, this map extends the classical Chern-Weil map $H^*_G(pt)\ra H^*_G(M)$.

\item{$\bullet$} If $M^G$ is nonempty, $\kappa_G:\H^*_G(\C)\ra \H^*_G(\cQ'(M))$ is injective. If $G$ is semisimple, $\H^*_G(\cQ'(M))$ is then a nontrivial conformal vertex algebra with Virasoro element $\kappa_G(\L)$ of central charge zero.

\item{$\bullet$} If the action of $G$ on $M$ is locally free, $\H^*_G(\cQ'(M)) _+ = 0$ and $\H^*_G(\cQ(M))_+ = 0$.

\item{$\bullet$} If $G$ is a torus $T$, $\H^*_T(\cQ'(M))_+ = 0$ if and only if the action of $T$ is locally free. The converse fails in general. For example, if $G$ is simple and $T\subset G$ is a torus, $\H^*_G(\cQ'(G/T))_+ = 0$ even though the action of $G$ on $G/T$ is not locally free. 

\item{$\bullet$} If $G$ is semisimple and $V$ is a faithful linear representation of $G$, $\H^*_G(\cQ(V))_+ = 0$.

\subsec{Outline of main results} In this paper, we continue the study of $\H^*_G(\cQ(M))$ and $\H^*_G(\cQ'(M))$. Our goal is to understand what kind of geometric information is contained in the positive-weight subspaces $$\H^*_G(\cQ(M))_+ = \bigoplus_{m>0} \H^*_G(\cQ(M))[m],~~~~~\H^*_G(\cQ'(M))_+ = \bigoplus_{m>0} \H^*_G(\cQ'(M))[m].$$ There are three basic results we need to establish. First, for $G$-invariant open sets $U,V\subset M$, there exist Mayer-Vietoris sequences
$$\cdots \ra \H^p_G(\cQ(U\cup V))\ra \H^p_G(\cQ(U))\oplus \H^p_G(\cQ(V))\ra \H^{p}_G(\cQ(U\cap V))\ra \cdots,$$ 
$$\cdots \ra \H^p_G(\cQ'(U\cup V))\ra \H^p_G(\cQ'(U))\oplus \H^p_G(\cQ'(V))\ra \H^{p}_G(\cQ'(U\cap V))\ra \cdots.$$ 

Second, $\H^*_G(\cQ'(-))$ is invariant under $G$-equivariant homotopy. That is, if $M$ and $N$ are $G$-manifolds and $\phi_0,\phi_1:M\ra N$ are equivariantly homotopic $G$-maps, the induced maps $\phi_0^*,\phi^*_1: \H^*_G(\cQ'(N))\ra\H^*_G(\cQ'(M))$ are the same.

Third, for any $G$ and $M$, $\H^*_G(\cQ'(M))$ and $\H^*_G(\cQ(M))$ have  {\it quasi-conformal structures}. That is, they admit an action of the subalgebra of the Virasoro algebra generated by $\{L_n|~n\geq -1\}$, such that $L_{-1}$ acts by $\partial$ and $L_0$ acts by $n\cdot  id$ on the subspace of weight $n$. The quasi-conformal structure provides a powerful vanishing criterion for $\H^*_G(\cQ'(M))_+$ and $\H^*_G(\cQ(M))_+$; it suffices to show that $L_0$ acts by zero.

Using these three basic tools, our goal will be to give a {\it relative} description of $\H^*_G(\cQ'(M))_+$ and $\H^*_G(\cQ(M))_+$ in terms of the vertex algebras $\H^*_K(\C)$ for connected normal subgroups $K$ of $G$, together with geometric data about $M$. If $K$ is abelian, Theorem 6.1 of \LLI~gives a complete description of $\H^*_K(\C)$, but if $K$ is non-abelian $\H^*_K(\C)$ is still a rather mysterious object. Computer calculations in the cases $K = SU(2)$ and $K = SU(3)$ indicate that $\H^*_K(\C)$ has a rich structure and contains many elements beyond the Virasoro element $\L$ that have no classical analogues. 

Since $G$-manifolds locally look like vector bundles over homogeneous spaces, a basic problem is to compute $\H^*_G(\cQ'(G/H))$ for any closed subgroup $H\subset G$.

\theorem{For any compact, connected $G$ and closed subgroup $H\subset G$,
\eqn\hspace{\H^*_G(\cQ'(G/H)) \cong \H^*_{K_0}(\C)\otimes H^*_{G'}(G/H),} where $K_0$ is the identity component of $K=Ker(G\ra Diff(G/H))$, and $G'=G/K_0$. Here $H^*_{G'}(G/H)$ is regarded as a vertex algebra in which all circle products are trivial except $\circ_{-1}$, and \hspace~is a vertex algebra isomorphism. }
\thmlab\homospace

Theorem \homospace~generalizes Corollary 6.14 of \LLSI, which deals with the case where $G$ is semisimple and $H$ is a torus. A consequence of this result is that for compact $M$, the degree $p$ and weight $n$ subspace $\H^p_G(\cQ'(M))[n]$ is finite-dimensional for all $p\in\Z$ and $n\geq 0$, which extends a well-known classical result in the case $n=0$. 

Next, we study $\H^*_G(\cQ(M))$ via the map $\H^*_G(\cQ'(M))\ra \H^*_G(\cQ(M))$ induced by the inclusion $\cQ'(M)\hookrightarrow\cQ(M)$. 

\theorem{For any $G$-manifold $M$,  
$$\H^*_G(\cQ(M)) \cong \H^*_{K_0}(\C)\otimes H^*_{G'}(M),$$ where $K_0$ is the identity component of $K = Ker(G\ra Diff(M))$ and $G' = G/K_0$. In particular, $\H^*_G(\cQ(M))_+ = 0$ whenever $K$ is finite.}
\thmlab\Qstructure

Thus $\H^*_G(\cQ(M))_+$ depends only on $K_0$, so it carries no other geometric information about $M$. An important consequence is that for any $G$ and $M$, $\H^*_G(\cQ(M))_+$ and $\H^*_G(\cQ'(M))_+$ are {\it vertex algebra ideals}, i.e., they are closed under $\alpha\circ_n$ and $\circ_n\alpha$ for all $n\in\Z$ and $\alpha$ in $\H^*_G(\cQ(M))$, $\H^*_G(\cQ'(M))$, respectively. 

Next, we study $\H^*_G(\cQ'(M))$, which in contrast to $\H^*_G(\cQ(M))$, carries non-trivial geometric information about $M$ beyond weight zero. We focus on three special cases: $G$ simple, $G = G_1\times G_2$ where $G_1,G_2$ are simple, and $G$ abelian. 

\theorem{(Positive-weight localization for simple group actions) For any simple $G$ and $G$-manifold $M$, the map \eqn\simp{ \H^*_G(\cQ'(M))_+\ra \H^*_G(\cQ'(M^G))_+} induced by the inclusion $i:M^G\ra M$, is an isomorphism of vertex algebra ideals. Hence $\H^*_G(\cQ'(M))_+ \cong \H^*_G(\C)_+\otimes H^*(M^G)$. Moreover, both the ring structure of $H^*(M^G)$ and the map $i^*:H^*_G(M)\ra H^*_G(M^G)$ are encoded in the vertex algebra structure of $\H^*_G(\cQ'(M))$.}
\thmlab\simpleloc

Using Theorem \simpleloc, together with results of R. Oliver \OI\OII~which describe the fixed-point subsets of group actions on contractible spaces, we prove our main theorem: 

\theorem{For any simple $G$, there is a sphere with infinitely many smooth actions of $G$, which have pairwise distinct chiral equivariant cohomology, but identical classical equivariant cohomology.}\thmlab\mainresult

One can even construct a morphism $f:M\ra N$ in the category of compact $G$-manifolds which induces a ring isomorphism $H^*_G(N)\ra H^*_G(M)$ (with $\Z$-coefficients), such that $\H^*_G(\cQ'(M))\neq \H^*_G(\cQ'(N))$. {\it Hence $\H^*_G(\cQ'(-))$ is a strictly stronger invariant than $H^*_G(-)$ on the category of compact $G$-manifolds.}

Similarly, in the case $G = G_1\times G_2$ where $G_1,G_2$ are simple, we describe $\H^*_G(\cQ'(M))$ in terms of the vertex algebras $\H^*_{G_1}(\C)$, $\H^*_{G_2}(\C)$ and the rings $H^*(M^G)$, $H^*_{G_1}(M^{G_2})$, and $H^*_{G_2}(M^{G_1})$.

The case where $G$ is the circle $S^1$ is analogous to the case of simple $G$. When $G$ is a general torus $T$, $\H^*_T(\cQ'(M))_+$ will typically depend on the family of rings $H^*_{T/T'}(M^{T'})$ for all subtori $T'\subset T$ for which $M^{T'}$ is non-empty, and can be quite complicated. As an example, we compute $\H^*_{T}(\cQ'(\C\P^2))$, where $T = S^1\times S^1$ and $\C\P^2$ is equipped with the usual linear action.

We conclude with a few remarks about $\H^*_G(\C)$. In \FF, it was suggested that the semi-infinite Weil complex $\cW(\gg)$ should play a role in semi-infinite geometry analogous to the role of the classical Weil complex $W(\gg)$. Note that $H^*_G(pt) = S(\gg^*)^G$ can be regarded either as the basic cohomology $H^*_{bas}(W(\gg))$, or as the Lie algebra cohomology $H^0(\gg,S(\gg^*))$. The analogue of $H^0(\gg,S(\gg^*))$ is the {\it semi-infinite cohomology} $H^{\infty + *}(\hat{\gg},\cS(\gg))$ \F\FGZ\Vor, whereas the analogue of $H^*_{bas}(W(\gg))$ is $\H^*_{bas}(\cW(\gg)) = \H^*_G(\C)$. Here $\cS(\gg)$ is the semi-infinite symmetric algebra on $\gg$. In contrast to the classical case, $\H^*_G(\C)\neq H^{\infty + *}(\hat{\gg},\cS(\gg))$. It would be interesting to construct an equivariant cohomology theory for manifolds in which $H^{\infty + *}(\hat{\gg},\cS(\gg))$ plays the role of $S(\gg^*)^G$ (as suggested in \FF), and compare it to our theory.

{\it Acknowledgement.} We thank G. Schwarz for helpful discussions on invariant theory and for pointing out to us Oliver's construction of group actions on disks for which the fixed-point sets have interesting homotopy types.

\subsec{General remarks about group actions on manifolds}

Let $G_0$ denote the identity component of $G$. For any $G$-manifold $M$, the finite group $\Gamma = G/G_0$ acts on the complex $(\cW(\gg)\otimes \cQ(M))_{G_0-bas}$ by DVA automorphisms, and we have
$$(\cW(\gg)\otimes \cQ(M))_{G-bas} = ((\cW(\gg)\otimes \cQ(M))_{G_0-bas})^{G} = ((\cW(\gg)\otimes \cQ(M))_{G_0-bas})^{\Gamma} .$$
Since the differential $d_{\cW} + d_{\cQ}$ commutes with the action of $\Gamma$ on $(\cW(\gg)\otimes \cQ(M))_{G_0-bas}$, $\Gamma$ acts on $\H^*_{G_0}(\cQ(M))$, and we have
$\H^*_G(\cQ(M)) = \H^*_{G_0}(\cQ(M))^{\Gamma}$. Similarly, $\H^*_G(\cQ'(M)) = \H^*_{G_0}(\cQ'(M))^{\Gamma}$. {\it Hence there is essentially no new content in studying $\H^*_G(-)$ for disconnected groups, and for the remainder of this paper, we will only consider the functor $\H^*_G(-)$ for connected $G$}.

We say that $G$ acts {\it effectively} on $M$ if $K = Ker(G\ra Diff(M))$ is trivial, and we say that $G$ acts {\it almost effectively} if $K$ is finite. Let $K_0$ denote the identity component of $K$ and let $G'=G/K_0$. 

\lemma{Let $M$ be a $G$-manifold, and suppose that $K = Ker(G\ra Diff(M))$ has positive dimension. Then $G'$ acts almost effectively on $M$ and 
$$\H^*_G(\cQ(M)) = \H^*_{K_0}(\C)\otimes \H^*_{G'}(\cQ(M)),~~~~ \H^*_G(\cQ'(M)) = \H^*_{K_0}(\C)\otimes \H^*_{G'}(\cQ'(M)).$$}\thmlab\effective

\proof Clearly $K/K_0 = Ker (G'\ra Diff(M))$, which is finite because $K$ is compact. Since $\gg = \gk\oplus \gg'$ where $\gk$ and $\gg'$ are the Lie algebras of $K$ and $G'$, respectively, we have $\cW(\gg) = \cW(\gk)\otimes \cW(\gg')$. Then 
$$(\cW(\gg)\otimes \cQ(M))_{G-bas} = \cW(\gk)_{K_0-bas}\otimes (\cW(\gg ')\otimes \cQ(M))_{G'-bas}.$$
Note that the differential $d = d_{\cW(\gg)} + d_{\cQ}$ of $\cW(\gg)\otimes \cQ(M)$ can be written as $d_{\cW(\gk)} + d_{\cW(\gg')} + d_{\cQ}$, and these three terms pairwise commute. Since $d_{\cW(\gk)}$ only acts on $\cW(\gk)_{K_0-bas}$ and $d_{\cW(\gg')} + d_{\cQ}$ only acts on $(\cW(\gg')\otimes \cQ(M))_{G'-bas}$, the claim follows. The proof of the corresponding statement for $\H^*_G(\cQ'(M))$ is the same. $\Box$

Any compact connected $G$ has a finite cover of the form $\tilde{G} = G_1\times \cdots\times G_k\times T$, where the $G_i$ are simple and $T$ is a torus. If $M$ is a $G$-manifold, the action can be lifted to $\tilde{G}$ so that $\Gamma = Ker(\tilde{G}\ra G)\subset Ker(\tilde{G}\ra Diff(M))$. 
\lemma{Suppose $M$ is a $G$-manifold, $K = Ker (G\ra Diff(M))$, and $\Gamma$ is a finite subgroup of $K$. Then $G/\Gamma$ acts on $M$ and 
$$\H^*_G(\cQ'(M)) = \H^*_{G/\Gamma}(\cQ'(M)),~~~~~ \H^*_G(\cQ(M)) = \H^*_{G/\Gamma}(\cQ(M)).$$
In particular, if $G$, $\tilde{G}$, and $\Gamma$ are as above, we have 
$$\H^*_G(\cQ'(M)) = \H^*_{\tilde{G}}(\cQ'(M)),~~~~~ \H^*_G(\cQ(M)) = \H^*_{\tilde{G}}(\cQ(M)).$$}
\thmlab\finitecover
\proof
$\Gamma$ acts trivially on $\cQ(M)$ since it acts trivially on $M$. The adjoint and coadjoint actions of $\Gamma$ on $\gg$ and $\gg^*$ are trivial, so $\Gamma$ also acts trivially on $\cW(\gg)$. Thus the $G$-invariance and the $G/\Gamma$-invariance conditions on $\cW(\gg)\otimes\cQ'(M)$ and $\cW(\gg)\otimes\cQ(M)$ are the same. Since $G$ and $G/\Gamma$ have the same Lie algebra, the $\gs\gg[t]$-basic condition is also the same. $\Box$

\newsec{Mayer-Vietoris Sequences}
In this section, we show that for $G$-invariant open sets $U,V\subset M$, there exist Mayer-Vietoris sequences 
$$\cdots\ra\H_G^p(\cQ'(U\cup V))\ra\H_G^p(\cQ'(U))\oplus\H_G^p(\cQ'(V))\ra\H_G^p(\cQ'(U\cap  
V))\ra\cdots,$$
$$\cdots\ra\H_G^p(\cQ(U\cup V))\ra\H_G^p(\cQ(U))\oplus\H_G^p(\cQ(V))\ra\H_G^p(\cQ(U\cap  
V))\ra\cdots.$$

\lemma{Let $M$ be a manifold and let $\gh$ be a Lie algebra. Suppose that $\cF$ is a sheaf of $C^\infty$-modules and of $\gh$-modules on $M$, where the two module structures are compatible, i.e. $\gh$ acts compatibly on the sheaf $C^\infty$ by derivations. Assume  
that $U,V$ are open sets, and that $\phi_U,\phi_V\in C^\infty(M)$ form  
an $\gh$-invariant partition of unity for the cover $\{U,V\}$ of  
$U\cup V$. If the invariant functor $(-)^\gh$ is applied to the standard exact  
sequence
$$
0\ra\cF(U\cup V)\ra\cF(U)\oplus\cF(V)\ra\cF(U\cap V)\ra 0
$$
the result is an exact sequence.}
\proof
Note that the left exactness of the standard sequence is just  the sheaf axiom, but the surjectivity of the last map is not true for a
general sheaf unless it is a fine sheaf (i.e. has the partition of  
unity property for sheaves.) Since $\cF$ is assumed to be a $C^\infty$  
sheaf, the exactness of our standard sequence is guaranteed.

Since the invariant functor $(-)^\gh$ is left exact, applying it to the  
standard sequence yields a left exact sequence. So it remains to show  
that
\eqn\invt{\cF(U)^\gh\oplus\cF(V)^\gh\ra\cF(U\cap V)^\gh}
is onto. Since the $C^\infty$ and $\gh$ structures on $\cF$ are  
compatible, the map is a $C^\infty$-module map, i.e. it is compatible  
with multiplications by functions. Let $a\in\cF(U\cap V)^\gh$ and let $\{\phi_U,\phi_V\}$ be an $\gh$-invariant partition of unity subordinate to the cover  
$\{U,V\}$ of $U\cup V$.

We claim that there is an extension of $(\phi_U|_{U\cap V})a\in\cF(U\cap  
V)^\gh$ by zero to all of $V$ (but not to all of $U$). We have
$$
(U\cap V)\cup (V\backslash supp(\phi_U))=V.
$$
For if $x\in V\backslash(U\cap V)$ then $x\notin supp(\phi_U)\subset  
U$, and so $x$ lies on the left side. Now to see that $(\phi_U|_{U\cap  
V})a\in\cF(U\cap V)^\gh$ and the zero $0\in\cF(V\backslash  
supp(\phi_U))^\gh$ glue together to form a section in $\cF(V)^\gh$, it  
is enough to check that $(\phi_U|_{U\cap V})a$ restricts to zero on the  
overlap $W=(U\cap V)\cap (V\backslash supp(\phi_U))$.
This restriction is equal to $(\phi_U|_{W})(a|_{W})$ since function  
multiplication commutes with restriction.
But $\phi_U|_{W}=0$ because $W\cap supp(\phi_U)=\emptyset$. This proves  
our claim. Call this extension $a_V\in\cF(V)^\gh$. Likewise let  
$a_U\in\cF(U)^\gh$ be the extension of $(\phi_V|_{U\cap V})a\in\cF(U\cap  
V)^\gh$ by zero to all of $U$.

Now under \invt~we have
$$
(a_V,-a_U)\mapsto (a_V|_{U\cap V})+(a_U|_{U\cap V})=(\phi_U|_{U\cap V})a+(\phi_V|_{U\cap V})a=
((\phi_U+\phi_V)|_{U\cap V})a=a
$$
since $\phi_U+\phi_V=1$. This proves that \invt~is onto. $\Box$

\remark{This result holds if we replace \lq\lq sheaf" with \lq\lq weak sheaf" since the reconstruction axiom  \recon~holds for finite covers.}

\theorem{Let $M$ be a $G$-manifold and let $U,V$ be $G$-invariant open  
sets in $M$. Then
$U\cap V\rightrightarrows U\coprod V\ra U\cup V$ induces a long exact sequence
$$
\cdots\ra\H_G^p(\cQ'(U\cup V))\ra\H_G^p(\cQ'(U))\oplus\H_G^p(\cQ'(V))\ra\H_G^p(\cQ'(U\cap  
V))\ra\cdots.
$$}
\proof
Regard $\cW=\cW(\gg)$ as a constant sheaf of vector spaces over $M$.  
Then $\cW\otimes\cQ'$ is a weak sheaf of $C^\infty$-modules where functions  
act only on the right factor. It is a weak sheaf of modules over the Lie  
algebra $\gs\gg[t]$, as shown in Section 3 of \LLSI. Choose a partition of  
unity $\phi_U,\phi_V$ of $U\cap V$ as before. Since $U,V$ are $G$-invariant sets, by averaging over $G$, we can assume that the two  
functions are $G$-invariant. Note that even though $\gs\gg[t]$ does  
not act by derivations on a general element of $\cW\otimes \cQ'$, $\gs\gg[t]$ does act by  
derivations on weight zero elements. Moreover, any  
$G$-invariant function $f$, regarded as $1\otimes f\in\cW\otimes\cQ'$, is chiral basic (i.e. $\gs\gg[t]$-invariant). So $\phi_U,\phi_V$ form an $\gs\gg[t]$-invariant partition of unity. Hence the preceding lemma can be applied to $\cF=\cW\otimes\cQ'$ and  
$\gh=\gs\gg[t]$. The invariant functor applied to the standard  
sequence for $\cF$ yields an exact sequence
$$
0\ra\cF(U\cup V)_{bas}\ra\cF(U)_{bas}\oplus\cF(V)_{bas}\ra\cF(U\cap  
V)_{bas}\ra 0,
$$ which induces the corresponding long exact sequence for chiral equivariant cohomology. $\Box$

The theorem holds if we replace $\cQ'$ by $\cQ$. The only tricky part  
here is that $\cQ$ is no longer a $C^\infty$-module because the Wick  
product on $\cQ$ is not associative. But we still have a $\C$-bilinear  
operation $C^\infty\times\cQ\ra\cQ,~(f,a)\mapsto :fg:$ which is a  
homomorphism of weak sheaves. Moreover, even though $\cQ$ is not functorial  
under general smooth mappings, it is functorial with respect to open  
inclusions. A partition of unity argument shows that the standard sequence for $\cQ$ is still exact. The proof of the preceding lemma then carries over to the case  
$\cF=\cW\otimes\cQ$.

\newsec{Homotopy Invariance of $\H^*_G(\cQ'(M))$}

Let $M$ and $N$ be $G$-manifolds, and let $\phi_0,\phi_1:M\ra N$ be $G$-equivariant maps. Let $I$ denote the interval $[0,1]$, which we regard as a $G$-manifold equipped with the trivial action. 

\definition{ A $G$-equivariant homotopy from $\phi_0$ to $\phi_1$ is a smooth $G$-equivariant map $\Phi:M\times I\ra N$ such that for all $x\in M$, $\Phi(x,0) = \phi_0(x)$ and $\Phi(x,1) = \phi_1(x)$. For each $t\in I$, $\phi_t:M\ra N$ will denote the map $\phi_t(x) = \Phi(x,t)$.} 

The main result in this section is
\theorem{Let $M$ and $N$ be $G$-manifolds, and let $\phi_0,\phi_1:M\ra N$ be $G$-equivariant maps. If there exists a $G$-equivariant homotopy $\Phi$ from $\phi_0$ to $\phi_1$ as above, the induced maps $\phi_0^*,\phi_1^*: \H^*_G(\cQ'(N))\ra \H^*_G(\cQ'(M))$ are the same.} 
\thmlab\hominv

We first define an appropriate notion of chiral chain homotopy in the category of $\gs\gg[t]$-modules, and show that two morphisms of $\gs\gg[t]$-modules which are chiral chain homotopic induce the same map in chiral equivariant cohomology. In the geometric setting, if $\Phi$ is a $G$-equivariant homotopy between $\phi_0,\phi_1:M\ra N$, we will construct a chiral chain homotopy between the induced maps $\phi_0^*,\phi_1^*:\cQ'(N)\ra\cQ'(M)$.

An immediate consequence of Theorem \hominv~is that if $M$ is $G$-equivariantly contractible to a submanifold $M'$, then $\H^*_G(\cQ'(M)) \cong \H^*_G(\cQ'(M'))$. By contrast, the functor $\H^*_G(\cQ(-))$ does {\it not} have this property. Let $G$ be simple and let $V$ be a faithful linear representation of $G$. Then $V$ is $G$-equivariantly contractible to the origin $o\in V$, but $\H^*_G(\cQ(V))_+ = 0$ whereas $\H^*_G(\cQ(o))_+ =\H^*_G(\C)_+\neq 0$.

\subsec{Chiral chain homotopies}
Suppose that $\cA$ and $\cB$ are $\gs\gg[t]$-modules. We define a {\it chiral chain homotopy} to be a linear map $P:\cA\ra \cB$, homogeneous of weight $0$ and degree $-1$, which is $G$-equivariant and satisfies
\eqn\chainhomo{P\iota_{\xi}^{\cA}(k) + \iota_{\xi}^{\cB}(k)P = 0,~~~~ P L_{\xi}^{\cA}(k) -L_{\xi}^{\cB}(k)P = 0,} for all $\xi\in\gg$ and $k\geq 0$.

\lemma{If $P:\cA\ra\cB$ is a chiral chain homotopy, the map $\tau= P d_{\cA}+d_{\cB} P$ is a morphism of $\gs\gg[t]$-modules.}

\proof The argument is similar to the proof of Proposition 2.3.1 in \GS. 
First note that $d_{\cB} \tau = d_{\cB} P d_{\cA} = \tau d_{\cA}$. From the assumption $L_{\xi}^{\cB}(k)P = P L_{\xi}^{\cA}(k)$, it is immediate that $$L_{\xi}^{\cB}(k) \tau = \tau L_{\xi}^{\cA}(k).$$ Finally, 
$$\iota_{\xi}^{\cB}(k)\tau = \iota_{\xi}^{\cB}(k) d_{\cB} P+\iota_{\xi}^{\cB}(k) Pd_{\cA} = \iota_{\xi}^{\cB}(k) d_{\cB} P - P \iota^{\cA}_{\xi}(k) d_{\cA} $$ $$= - d_{\cB}\iota_{\xi}^{\cB}(k) P + L_{\xi}^{\cB}(k) P + P d_{\cA} \iota_{\xi}^{\cA}(k) - P L_{\xi}^{\cA}(k) = (d_{\cB} P + P d_{\cA})\iota_{\xi}^{\cA}(k) = \tau \iota_{\xi}^{\cA}(k). \Box $$

We say that two $\gs\gg[t]$-module homomorphisms $\phi_0,\phi_1:\cA\ra \cB$ are {\it chiral chain homotopic} if there is a chiral chain homotopy $P:\cA\ra \cB$ such that $\phi_1-\phi_0 = \tau$. This clearly implies that $\phi_0,\phi_1$ induce the same map from $\H^*_{bas}(\cA)\ra \H^*_{bas}(\cB)$.

\lemma{Chiral chain homotopic maps $\phi_0,\phi_1:\cA\ra \cB$ induce the same map from $\H^*_G(\cA)\ra \H^*_G(\cB)$.}

\proof This is analogous to Proposition 2.4.1 of \GS. The map $id\otimes P: \cW\otimes\cA \ra \cW\otimes\cB$ is a chiral chain homotopy between $id\otimes \phi_0$, $id \otimes \phi_1$. Hence $id\otimes\phi_0$, $id\otimes \phi_1$ induce the same map from $\H^*_G(\cA) = \H^*_{bas}(\cW\otimes\cA)\ra  \H^*_{bas}(\cW\otimes\cB) = \H^*_G(\cB)$. $\Box$

Suppose that $\phi_0, \phi_1:M\ra N$ are $G$-equivariantly homotopic via $\Phi:M\times I\ra N$. We recall the classical construction of a chain homotopy $P:\Omega^*(N)\ra \Omega^{*-1}(M)$ between $\phi_0^*,\phi_1^*:\Omega^*(N)\ra \Omega^*(M)$, following \GS. For fixed $x\in M$, consider the curve in $N$ given by $s\mapsto \phi_s(x)$, and let
$\xi_t:M\ra TN$ be the map which assigns to $x$ the tangent vector to this curve at $s=t$. 
Consider the map \eqn\deff{f_t:\Omega^*(N)\ra \Omega^{*-1}(M),~~~~~\sigma\mapsto \phi_t^*(\iota_{\xi_t}(\sigma)).} At each $x\in M$, this map is defined by
\eqn\phiiota{\phi_t^*(\iota_{\xi_t}(\sigma))(\eta_1,\dots,\eta_k) = \sigma(\xi_t(x),d\phi_t(\eta_1),\dots,d\phi_t(\eta_k)),} given vectors $\eta_1,\dots,\eta_k\in TM_x$.

A well-known formula (see \GS) asserts that
\eqn\diffcalc{\frac{d}{dt} \phi_t^*\sigma = \phi_t^*(\iota_{\xi_t} (d\sigma) + d (\phi_t^*(\iota_{\xi_t}(\sigma)).}
Define $P: \Omega^*(N)\ra \Omega^{*-1}(M)$ by $P\sigma = \int_0^1 \phi_t^*(\iota_{\xi_t}(\sigma))~dt$. Integrating \diffcalc~over $I$ shows that $P d + d P = \phi_1^*-\phi_0^*$, so $P$ is a chain homotopy. Since $\phi_0,\phi_1,\Phi$ are $G$-equivariant maps, $f_t$ is $G$-equivariant and satisfies $ f_t \iota^N_{\xi} +\iota^M_{\xi} f_t = 0$, for all $t\in I$. It follows that $P$ is also $G$-equivariant and satisfies $P\iota^N_{\xi} + \iota^M_{\xi} P = 0$. Hence $P$ is a chain homotopy in the sense of \GS, and $\phi_0,\phi_1$ induce the same maps in equivariant cohomology.

We need to show that $P$ extends to a map $\P: \cQ'(N)\ra \cQ'(M)$ which is a chiral chain homotopy between $\phi_0^*,\phi_1^*:\cQ'(N)\ra\cQ'(M)$. By Lemma 3.2 of \LLSI, for each $m\geq 0$, we may regard $\cQ'_M[m]$ as a smooth vector bundle over $M$ of finite rank, which has a local trivialization induced from a collection of charts on $M$. Given a coordinate open set $U\subset M$ with coordinates $\gamma^1,\dots,\gamma^n$, $U\times V$ is a local trivialization of $\cQ'_M$, where $V$ is the vector space with basis consisting of all nonzero monomials of the form
\eqn\monomial{\partial^{r_1}\gamma^{i_1}\cdots\partial^{r_k}\gamma^{i_k}\partial^{s_1} c^{j_1}\cdots \partial^{s_l} c^{j_l},} where $ r_1+\cdots + r_k+s_1+\cdots + s_l =m$, and each $r_i>0$ and $ s_i\geq 0$. Here $c^i$ denotes the coordinate one-form $d_{\cQ}\gamma^i$.

Let $\pi:M\times I \ra M$ be the projection onto the first factor. We can pull back $\cQ'_M[m]$ to a vector bundle $\pi^*(\cQ'_M[m])\ra M\times I$. Let $\Gamma[m] = \Gamma(M\times I, \pi^*(\cQ'_M[m]))$ denote the space of smooth sections $$\sigma:M\times I\ra \cQ'_M[m],~~~~\sigma(x,t) = (x,v(x,t)),$$ where in local coordinates $v(x,t) = \sum_{j\in J} f_j(\gamma^1,\dots,\gamma^n,t) \mu_j$. Here the set $J$ indexes all monomials of the form \monomial, and each $f_j$ is a smooth function on $U\times I$. 

Note that $\Gamma= \oplus_{m\geq 0} \Gamma[m]$ is an $\gs\gg[t]$-algebra, and that $\frac{d}{dt}$ and $\partial$ are derivations on $\Gamma$. It is clear from the local description of $\sigma(x,t)$ that $$\frac{d}{dt} \partial (\sigma(x,t)) = \partial (\frac{d}{dt}\sigma(x,t)).$$ Furthermore, the (fiberwise) integral $\int_0^1 \sigma(x,t) dt$ is a well-defined map from $\Gamma\ra \Gamma(M,\cQ'_M)$, and $\int_0^1 \partial \sigma(x,t) dt = \partial \int_0^1 \sigma(x,t) dt$.

\subsec{Relative derivations}
Let $\phi:A\ra B$ be a map of $G^*$-algebras. A degree-homogeneous map $f:A\ra B$ is called a {\it derivation relative to $\phi$}, or a $\phi$-derivation if \eqn\rel{f(ab) = f(a)\phi(b) + (-1)^{(deg~f)( deg~a)}\phi(a)f(b),} for all homogeneous $a,b\in A$. For example, in the case $A = \Omega(N)$, $B = \Omega(M)$, $\phi = \phi^*_t$, the map $f_t$ given by \deff~is a $\phi^*_t$-derivation.

Similarly, given a morphism $\phi:\cA\ra \cB$ of $\gs\gg[t]$-algebras, a degree-weight homogeneous linear map $f:\cA\ra \cB$ will be called a $\phi$-derivation if 
\eqn\relder{f(a\circ_n b) = f(a)\circ_n \phi(b) + (-1)^{(deg~f )(deg~a)}\phi(a)\circ_n f(b),} for all homogeneous $a,b\in \cA$ and $n\in\Z$. Clearly $f(1) = 0$ and $f (\partial a) = \partial f(a)$ for all $a\in \cA$. If $\cA,\cB$ are abelian vertex algebras, to check that $f$ is a $\phi$-derivation, it is enough to show that for all $a,b\in\cA$,
\eqn\abcase{f(:ab:) = ~:f(a)\phi(b):~ + (-1)^{(deg~f)( deg~a)}: \phi(a) f(b):~,~~~~~~~ f(\partial a) = \partial f(a).} 

\remark{A $\phi$-derivation $f$ is determined by its values on a set of generators of $\cA$. In the case $\cA = \cQ'(N)$, $\cB = \cQ'(M)$, $\phi = \phi^*_t$, since $\cQ'(N)$ is generated by $\Omega(N)$, any two $\phi^*_t$-derivations which agree on $\Omega(N)$ must agree on all of $\cQ'(N)$.}
\thmlab\remarkI

\remark{Suppose that $f$ is a $\phi$-derivation and $\delta_{\cA},\delta_{\cB}$ are vertex algebra derivations on $\cA,\cB$, respectively, which are homogeneous of degree $d$ and satisfy $\phi\circ \delta_{\cA} = \delta_{\cB}\circ \phi$. Then $f\circ \delta_{\cA} -(-1)^{(deg~ f)(d)} \delta_{\cB}\circ f$ is also a $\phi$-derivation. }
\thmlab\remarkII

For example, for any $\xi\in\gg$, the operators $\iota_{\xi}^{\cA}(0)$, $\iota_{\xi}^{\cB}(0)$ are vertex algebra derivations of degree $-1$, and $L_{\xi}^{\cA}(0)$, $L_{\xi}^{\cB}(0)$ are vertex algebra derivations of degree $0$. Hence $$f\circ \iota_{\xi}^{\cA}(0) -(-1)^{(deg ~f)} \iota_{\xi}^{\cB}(0)\circ f,~~~~~~~~~f\circ L_{\xi}^{\cA}(0)- L_{\xi}^{\cB}(0)\circ f$$ are $\phi_t^*$-derivations. If $\cA$ and $\cB$ are abelian vertex algebras, the operators $\iota_{\xi}^{\cA}(k)$, $\iota_{\xi}^{\cB}(k)$, $L_{\xi}^{\cA}(k)$, $L_{\xi}^{\cB}(k)$ are vertex algebra derivations for all $k\geq 0$, so $$f\circ \iota_{\xi}^{\cA}(k) - (-1)^{(deg ~f)} \iota_{\xi}^{\cB}(k)\circ f, ~~~~~~~~~f\circ L_{\xi}^{\cA}(k) -L_{\xi}^{\cB}(k)\circ f$$ are $\phi^*_t$-derivations.

\lemma{There is a unique extension of the map $f_t:\Omega^*(N)\ra \Omega^{*-1}(M)$ defined by \deff~to a linear map $F_t: \cQ'^*(N)\ra \cQ'^{*-1}(M)$, which is a $\phi^*_t$-derivation.}
\thmlab\extension

\proof We first construct $F_t$ locally. On a coordinate open set $U\subset N$ with coordinates $\gamma^1,\dots,\gamma^n$, recall from \LLSI~that $\cQ'(U)$ is the abelian vertex algebra with generators $f\in \cC^{\infty}(U)$ and $c^i$, $i=1,\dots,n$, subject to the relations:
\eqn\relations{1(z) -1,~~~~~~~~~(fg)(z) - f(z)g(z),~~~~~~~~~ \partial f(z) - \frac{\partial f}{\partial \gamma^i}(z) \partial \gamma^i(z),} for all $f,g\in \cC^{\infty}(U)$. We define $F_t$ on generators by $F_t(f(\gamma^1,\dots,\gamma^n)) = 0$ and $F_t(c^i) = f_t (c^i)$,
and then extend $F_t$ to a linear map on all of $\cQ'(U)$ using the $\phi^*_t$-derivation property \abcase. Since the relations \relations~are all homogeneous of degree 0 and $F_t$ lowers degree by one, it is clear that $F_t$ annihilates these relations, and hence is well-defined. Using the fact that $F_t(f(\gamma^1,\dots,\gamma^n)) = 0$ for any $f\in \cC^{\infty}(U)$, it is easy to check that the definition of $F_t$ is coordinate-independent. 

Finally, cover $N$ with coordinate open sets $\{U_i\}$, and define $F_t|_{U_i}$ as above. Fix a partition of unity $\{\psi_i\}$ subordinate to this covering. For each $m\geq 0$, we define $F_t:\cQ'(N)[m]\ra \cQ'(M)[m]$ by $F_t = \sum_i \psi_i F_t|_{U_i}$, which is well-defined since $\cQ'(N)[m]$ is a fine sheaf. Moreover, $F_t$ still satisfies \abcase~on each $\cQ'(N)[m]$ because $F_t(\psi_i) = 0$. Finally, since $\cQ'(N) = \bigoplus_{m\geq 0} \cQ'(N)[m]$, this defines $F_t$ on all of $\cQ'(N)$. $\Box$

\lemma{ Equation \diffcalc~holds for any $\sigma$ in $\cQ'(N)$, not just $\Omega(N)$.}

\proof By the preceding lemma, both sides of \diffcalc~are well-defined. It suffices to show that it holds at $t=s$ for each $s\in I$. Let $$g  = \frac{d}{dt} \phi^*_t|_{t=s}( - ),~~~~~~ h = \phi_s^* (\iota_{\xi_s} (0) d(-)) + d\phi_s^*( \iota_{\xi_s}(0)(-)),$$ which are the maps from $\cQ'(N)\ra \cQ'(M)$ appearing on the left and right sides of \diffcalc, evaluated at $t=s$. Clearly $g$ and $h$ are both $\phi^*_s$-derivations, and since $g$ and $h$ agree on generators, they agree on all of $\cQ'(N)$ by Remark \remarkI. $\Box$

We now define $\P:\cQ'(N)\ra \cQ'(M)$ by $\P(\sigma) = \int_0^1 F_t(\sigma) dt$, which coincides with $P$ at weight zero. Integration of \diffcalc~over $I$ shows that $d\P+\P d = \phi_1^*-\phi_0^*$. Finally, we need to show that the map $\P$ constructed above is in fact a chiral chain homotopy. Recall that for all $\sigma\in\Omega(N)$, $\xi\in\gg$, and $t\in I$, $f_t$ satisfies
\eqn\commutator{ f_t L_{\xi}^N - L_{\xi}^M f_t = 0,~~~~~~~~~~~  f_t \iota_{\xi}^N + \iota_{\xi}^M f_t = 0.} 
For $\xi\in\gg$ and $k\geq 0$, consider the maps $$R_{t,\xi,k} = F_t L^N_{\xi}(k) - L^M_{\xi}(k) F_t, ~~~~~ S_{t,\xi,k}= F_t \iota^N_{\xi}(k) + \iota^M_{\xi}(k) F_t.$$ By Remark \remarkII, $R_{t,\xi,k}$ and $S_{t,\xi,k}$ are $\phi^*_t$-derivations from $\cQ'(N)\ra \cQ'(M)$, which are homogeneous of weight $-k$ and degree $-1$ and $-2$, respectively. For $k>0$, 
$R_{t,\xi,k}$ and $S_{t,\xi,k}$ both act by zero on $\cQ'(N)[0]$ by weight considerations. For $k=0$, $R_{t,\xi,k}$ and $S_{t,\xi,k}$ act by zero on $\cQ'(N)[0]$ by \commutator. Since $R_{t,\xi,k}$ and $S_{t,\xi,k}$ are $\phi^*_t$-derivations, it follows from Remark \remarkI~that they act by zero on all of $\cQ'(N)$.

Finally, since this holds for each $t\in I$, it is immediate that 
$$\P L_{\xi}^N(k) (\sigma)-L_{\xi}^M(k)\P(\sigma) = \int_0^1 R_{t,\xi,k}(\sigma) dt = 0,~~~~\P\iota_{\xi}^N(k)(\sigma) + \iota_{\xi}^M(k) \P(\sigma) = \int_0^1 S_{t,\xi,k}(\sigma) dt = 0,$$ for all $\xi\in\gg$, $k\geq 0$ and $\sigma\in\cQ'(N)$. Hence $\P$ is a chiral chain homotopy, as desired.

\newsec{A Quasi-conformal Structure on $\H^*_G(\cQ'(M))$ and $\H^*_G(\cQ(M))$}

When $G$ is semisimple, $\H^*_G(\C)$ is a conformal vertex algebra with Virasoro element $\L$ of central charge 0. For any $M$, $\H^*_G(\cQ'(M))$ is then a conformal vertex algebra with Virasoro element $\kappa_G(\L)$. The vanishing of $\kappa_G(\L)$ is a necessary and sufficient condition for the vanishing of $\H^*_G(\cQ'(M))_+$ since $\kappa_G(\L)\circ_1 \omega = n \omega$ for all $\omega\in\H^*_G(\cQ'(M))[n]$. Unfortunately, $\H^*_G(\cQ'(M))$ is not a conformal vertex algebra when $G$ has a positive-dimensional center, so a priori we have no such vanishing criterion for general $G$.

In this section we show that for any $G$ and $M$, both $\H^*_G(\cQ'(M))$ and $\H^*_G(\cQ(M))$ have a {\it quasi-conformal} structure, that is, an action of the subalgebra of the Virasoro algebra generated by $\{L_n|\ n\geq -1\}$, such that $L_{-1} = L\circ_0$ acts by $\partial$, and $L_0 = L\circ_1$ acts by $n\cdot id$ on the subspace of weight $n$. Thus we have a similar vanishing criterion for $\H^*_G(\cQ'(M))_+$ and $\H^*_G(\cQ(M))_+$ for any $G$; it suffices to show that $L\circ_1$ acts by zero.

We work in the setting of a general $O(\gs\gg)$ topological vertex algebra (TVA), which we defined in \LLSI. Recall that an $O(\gs\gg)$ TVA is a degree-weight graded DVA  $(\cA,d)$ equipped with an $O(\gs\gg)$-structure $(\xi,\eta)\mapsto L^{\cA}_\xi+\iota^{\cA}_\eta$, a chiral horizontal element $g^{\cA}$, such that $L^{\cA}=d g^{\cA}$ is a conformal structure, with respect to which $L^{\cA}_{\xi}$ and $\iota^{\cA}_{\eta}$ are primary of weight one. We call $g^{\cA}$ a chiral contracting homotopy of $\cA$. Given an $O(\gs\gg)$ TVA $(\cA,d)$, a differential vertex subalgebra $\cB$ is called a half $O(\gs\gg)$ TVA if the nonnegative Fourier modes of the vertex operators $\iota^{\cA}_\xi$ and $g^{\cA}$ preserve $\cB$. Note that the nonnegative Fourier modes of $L^{\cA}_\xi=d\iota^{\cA}_\xi$ and $L^{\cA}=d g^{\cA}$ automatically preserve $\cB$ as well. In particular, the action of $\{L^{\cA}\circ_n|\ n\geq 0\}$ is a quasi-conformal structure on $\cB$. Since $[d,g^{\cA}\circ_1] = L^{\cA}\circ_1$ and $g^{\cA}\circ_1$ acts on $\cB_{bas}$, Theorem 4.8 of \LLSI~shows that $\H^*_{bas}(\cB)$ vanishes beyond weight zero. 

For a $G$-manifold $M$, $\cQ(M)$ is our main example of an $O(\gs\gg)$ TVA. In local coordinates, 
\eqn\tvamanifold{g=g^M = b^i\partial \gamma^i,~~~~~L=L^M = \beta^i\partial\gamma^i - b^i\partial c^i.} The subalgebra $\cQ'(M)$ is then a half $O(\gs\gg)$ TVA as above. Recall that the semi-infinite Weil algebra $\cW = \cW(\gg)$ is {\it not} an $O(\gs\gg)$ TVA since there is no chiral horizontal element $g^{\cW}$ satisfying $d g^{\cW} = L^{\cW}$. 

Let $\cB$ be a half $O(\gs\gg)$ TVA inside some $O(\gs\gg)$ TVA $\cA$ as above. Then the non-negative Fourier modes of 
$$L^{tot} = L^{\cW}\otimes 1 + 1\otimes L^{\cA} \in \cW\otimes \cA$$ act on $\cW\otimes \cB$, giving $\cW\otimes \cB$ a quasi-conformal structure. Moreover, 
$$L^{tot}_{\xi} = L^{\cW}_{\xi}\otimes 1 + 1\otimes L^{\cA}_{\xi},~~~~~~\iota^{tot}_{\xi} = \iota^{\cW}_{\xi}\otimes 1 + 1\otimes \iota^{\cA}_{\xi}$$ are primary of weight one with respect to $L^{tot}$, and $d L^{tot} = 0$. 

\theorem{$L^{tot}\circ_n$ operates on $\H^*_G(\cB)$ for $n\geq 0$, and gives $\H^*_G(\cB)$ a quasi-conformal structure.  If $G$ is semisimple, $\kappa_G(\L)\circ_n= L^{tot}\circ_n $ as operators on $\H^*_G(\cB)$ for all $n\geq 0$.}\thmlab\qconf
\proof This is immediate from Theorem 4.8 and Theorem 4.17 of \LLSI. $\Box$

\subsec{A vanishing criterion}
For any compact $G$, and any $\cA$ and $\cB$ as above, $L^{tot}\circ_1$ acts by $n\cdot id$ on the subspaces $\cA[n]$ and $\cB[n]$ of weight $n$. Hence the vanishing of $L^{tot}\circ_1$ on $\H^*_G(\cB)$ (resp. $\H^*_G(\cA))$ is equivalent to the vanishing of $\H^*_G(\cB)_+$ (resp. $\H^*_G(\cA)_+$). The next lemma gives a useful vanishing criterion for $L^{tot}\circ_1$.

\lemma{Suppose that $\alpha\in \cW\otimes \cB$ is homogeneous of weight $2$ and degree $-1$, is $G$-invariant, chiral horizontal, and satisfies 
\eqn\circleone{L^{tot}_{\xi}\circ_1\alpha = \beta^{\xi}\otimes 1,} for all $\xi\in\gg$. Then $L^{tot}\circ_1$ acts by zero on $\H^*_G(\cB)$, and we have $\H^*_G(\cB)_+ = 0$.}
\thmlab\vanishing
\proof Recall from \LLI~that $d(\beta^{\xi_i} \partial c^{\xi'_i}\otimes 1) = L^{\cW}\otimes 1$, but $\beta^{\xi_i}\partial c^{\xi'_i}\otimes 1$ is not chiral horizontal since $\iota^{\cW}_{\xi} \circ_n (\beta^{\xi_i} \partial c^{\xi'_i}\otimes 1) = -\delta_{n,1} \beta^{\xi}\otimes 1$ for $n\geq 0$. Let $\omega_0 = \beta^{\xi_i} \partial c^{\xi'_i}\otimes 1 + d\alpha$. Clearly $d\omega_0= L^{\cW}\otimes 1$ and  
$$\iota^{tot}_{\xi}\circ_n (d\alpha) = L^{tot}_{\xi}\circ_n\alpha = \delta_{n,1} \beta^{\xi}\otimes 1 = - \iota^{tot}_{\xi}\circ_n (\beta^{\xi_i} \partial c^{\xi'_i}\otimes 1)$$ for $n\geq 0$, since $\alpha$ is chiral horizontal. It follows that $\iota^{tot}_{\xi}\circ_n \omega_0 = 0$ for all $n\geq 0$, so $\omega_0$ is chiral horizontal. In particular, the non-negative Fourier modes of $\omega_0$ act on $(\cW\otimes \cB)_{bas}$. Finally, let $\omega = \omega_0 + g^{\cA} \in \cW\otimes \cA$. Since $d g^{\cA} = L^{\cA}$ we have $d \omega = L^{tot}$. The non-negative Fourier modes of $\omega$ clearly preserve $(\cW\otimes \cB)_{bas}$ since both $\omega_0$ and $g^{\cA}$ have this property. In particular, $[d,\omega\circ_1] = L^{tot}\circ_1$, so $\omega\circ_1$ is a contracting homotopy for $L^{tot}\circ_1$, as desired. $\Box$ 

This lemma clearly holds if we replace $\cB$ with $\cA$. Note that when $G$ is semisimple, the existence of $\alpha$ as above guarantees that $\kappa_G(\L) = 0$; take $\omega = \beta^{\xi_i}\partial c^{\xi'_i} \otimes 1+ d\alpha + \theta^{\xi_i}_{\cS} b^{\xi_i}\otimes 1$. An OPE calculation shows that $\omega$ is chiral basic and $d\omega = \kappa_G(\L)$.

In \LLSI, we considered two situations where we can construct $\alpha$ as above. First, if $G$ acts locally freely on $M$, we have a map $\gg^*\ra \Omega^1(M)$ sending $\xi'\ra \theta^{\xi'}$, such that $\iota_{\xi}\theta^{\eta'} = \bra\eta',\xi\ket$. The $\theta^{\xi'}$ are known as {\it connection one-forms}. Choose an orthonormal basis $\{\xi_i\}$ for $\gg$ relative to the Killing form, and let $\Gamma^{\xi'_i} = g\circ_0 \theta^{\xi'_i}$. Then $$\alpha = \beta^{\xi_i}\otimes \Gamma^{\xi'_i}\in\cW\otimes \cQ'(M)$$ is $G$-invariant, chiral horizontal, and satisfies \circleone. This shows that $\H^*_G(\cQ'(M))_+ = 0$ and $\H^*_G(\cQ(M))_+ = 0$. 

Second, let $V$ be a faithful linear representation of $G$, and let $\rho:\gg\ra End(V)$ denote the corresponding representation of $\gg$. The induced bilinear form $\bra\xi,\eta\ket = Tr(\rho(\xi)\rho(\eta))$ is nondegenerate, so we may identify $\gg$ with $\gg^*$ via $\bra,\ket$ and fix an orthonormal basis $\xi_i$ of $\gg$.  Let $x_k$ be a basis of $V$ and $x'_k$ the corresponding dual basis for $V^*$. Define
\eqn\linearspace{\alpha = \beta^{\xi_i}\otimes \Gamma^{\xi_i} - \beta^{\xi_i} b^{\xi_j}\otimes \iota_{\xi_j}\circ_0 \Gamma^{\xi_i}\in\cW\otimes\cQ(V),}
where $\Gamma^{\xi_i} = \beta^{\rho(\xi_i)(x_k)}\gamma^{x'_k}$. An OPE calculation shows that $\alpha$ satisfies the conditions of Lemma \vanishing, so $\H^*_G(\cQ(V))_+ = 0$.

The next lemma shows that locally defined vertex operators $\alpha$ satisfying the conditions of Lemma \vanishing~can be glued together.

\lemma{Let $M$ be a $G$-manifold and let $\{U_{i}|\ i\in I\}$ be a cover of $M$ by $G$-invariant open sets. Suppose that $\alpha_i\in\cW\otimes \cQ'(U_i)$ satisfies the conditions of Lemma \vanishing. Then $\H^*_G(\cQ'(M))_+ = 0$.} \thmlab\locI

\proof Let $\{\phi_i|\ i\in I\}$ be a $G$-invariant partition of unity subordinate to the cover. Let $\alpha = \sum_i \phi_i \alpha_i$, which is a well-defined global section of $ \cW\otimes \cQ'(M)$. Moreover, since $\phi_i$ is basic, it follows that $\alpha$ remains $G$-invariant, $G$-chiral horizontal, and satisfies \circleone. $\Box$

\remark{Similarly, if $\alpha_i\in \cW\otimes \cQ(U_i)$ satisfies the conditions of Lemma \vanishing, $\alpha = \sum_i \phi_i \alpha_i\in \cW\otimes \cQ(M)$ does as well, so that $\H^*_G(\cQ(M))_+ =0$.}\thmlab\locII

\newsec{$\H^*_G(\cQ'(G/H))$ for Homogeneous Spaces $G/H$}
A basic fact about $G$-manifolds which can be found in \TD~is that locally they look like vector bundles over homogeneous spaces $G/H$. 

\theorem{Let $G$ be a compact Lie group and let $M$ be a smooth $G$-manifold. For each point  $x\in M$, the isotropy group $G_x$ is a closed subgroup of $G$ and the orbit $Gx$ is $G$-diffeomorphic to $G/G_x$. Moreover, $Gx$ has a $G$-invariant tubular neighborhood $U_x$ which is $G$-diffeomorphic to the bundle $G\times_{G_x} V$ for some real $G_x$-representation $V$.}
\thmlab\localstructure

By homotopy invariance, $\H^*_G(\cQ'(U_x)) \cong  \H^*_G(\cQ'(G/G_x))$. Thus the problem of computing $\H^*_G(\cQ'(G/H))$ for any closed subgroup $H\subset G$ is an important building block in the study of $\H^*_G(\cQ'(M))$ for general $M$.

Suppose first that $K = Ker (G\ra Diff(G/H))$ has positive dimension. Clearly $K\subset H$; as in Lemma \effective, let $K_0$ denote the identity component of $K$, and let $G' = G/K_0$ and $H' = H/K_0$. By Lemma \effective, we have $$\H^*_G(\cQ'(G/H)) = \H^*_{K_0}(\C)\otimes \H^*_{G'}(\cQ'(G/H)) = \H^*_{K_0}(\C)\otimes \H^*_{G'} (\cQ'(G'/H')).$$ Since $G'$ acts almost effectively on $G'/H'$, which is a homogeneous space for $G'$, we may assume without loss of generality that $K$ is finite. In this case, $\H^*_{K_0}(\C) = \C$, so in order to prove Theorem \homospace, it suffices to prove that $\H^*_G(\cQ'(G/H))_+ = 0$.

We need a basic property of any simple finite-dimensional complex Lie algebra $\gg$. Suppose that $\gh\subset \gg$ is a Lie subalgebra of positive codimension.  Via the Killing form, we identify $\gg$ with $\gg^*$, and in particular we identify $\gh^{\perp}$ with $(\gg/\gh)^*$. Note that $(\gg/\gh)^*$ is a representation of $\gh$. Regarding $(\gg/\gh)^*$ as a subspace of $\gg^*$, we may consider the subspace $ad^*(\gg/\gh)^*\subset\gg^*$ under the coadjoint action of $\gg$ on $\gg^*$.

\lemma{$ad^*(\gg/\gh)^* = \{ad_{\xi}^*(\eta')|\ \xi\in\gg,\eta'\in (\gg/\gh)^*\} = \gg^*$.}
\thmlab\surjectivity

\proof If $\gh\subset\gh'\subset\gg$ for some other subalgebra $\gh'$, we have $(\gg/\gh')^* \subset(\gg/\gh)^*$. Then $ad^*(\gg/\gh')^* \subset ad^*(\gg/\gh)^*$, so we may assume that $\gh$ is a maximal subalgebra of $\gg$ for which $\gh\neq \gg$. If $ad^*(\gg/\gh)^*\neq\gg^*$, there is some nonzero $\xi_0\in\gg$ such that
\eqn\adI{
\bra \xi_0, ad_{\xi_1}^*\eta'\ket = 0,\ \ \forall \xi_1\in\gg,\ \ \eta'\in(\gg/\gh)^*}

Let $B$ denote the set of $\xi_0$ which satisfy \adI. Then if $\xi_0\in B$, for any $\xi_1\in\gg$ and $\eta'\in(\gg/\gh)^*$ we have 

\eqn\adII{ \bra ad_{\xi_1}\xi_0,\eta'\ket = 0  \Longleftrightarrow  ad_{\xi_1}\xi_0\in\gh. }

Suppose first that $\xi_0\notin\gh$. Then $\gh' = \C\xi_0\oplus\gh$ is a Lie subalgebra of $\gg$. Since $\gh$ is maximal, $\gh' = \gg$. It follows from \adII~that $\gh$ is a nontrivial ideal of $\gg$, which contradicts the assumption that $\gg$ is simple. Hence $\xi_0\in\gh$, so we have $B\subset \gh$.

For any $\xi_1\in\gg$, we have the decomposition $\xi_1 = \xi_1^{\gh} + \xi_1^{\gh^{\perp}}$, where $\xi_1^{\gh}\in\gh$ and $\xi_1^{\gh^{\perp}}\in\gh^{\perp}$. Note that $ad_{\xi_1^{\gh^{\perp}}}\xi_0 = -ad_{\xi_0} \xi_1^{\gh^{\perp}}\in\gh^{\perp}$ since $\gh^{\perp} = (\gg/\gh)^*$ is a representation of $\gh$.
But $ad_{\xi_1^{\gh^{\perp}}}\xi_0\in\gh$ by \adII, so that $ad_{\xi_1^{\gh^{\perp}}}\xi_0\in\gh\cap \gh^{\perp}$ and we have $ad_{\xi_1^{\gh^{\perp}}}\xi_0 = 0$. Hence
 $$ad_{\xi_1}\xi_0 = ad_{\xi_1^{\gh}}\xi_0+ad_{\xi_1^{\gh^{\perp}}}\xi_0 = ad_{\xi_1^{\gh}}\xi_0.$$ 
 
 We claim that for any $\xi_2\in\gg$ and $\eta'\in(\gg/\gh)^*$, $$\bra ad_{\xi_1}\xi_0,ad^*_{\xi_2}\eta'\ket = 0,$$ so by definition $ad_{\xi_1}\xi_0\in B$. Hence $B$ is a nontrivial ideal since $B\neq 0$ and $B\subset\gh$, which is impossible.
 
 We have $$\bra ad_{\xi_1^{\gh}}\xi_0,ad^*_{\xi_2}\eta'\ket = - \bra \xi_0, ad^*_{[\xi_1^{\gh},\xi_2]}\eta'\ket - \bra \xi_0, ad^*_{\xi_2} ad^*_{\xi_1^{\gh}}\eta'\ket.$$
 
 The first term above is zero since $\xi_0\in B$. The second term is also zero since $ad^*_{\xi_1^{\gh}}\eta'\in(\gg/\gh)^*$, since $(\gg/\gh)^*$ is a representation of $\gh$. $\Box$
 
 \remark{Lemma \surjectivity~remains true if $\gg$ is semisimple and $\gh$ does not contain any simple component of $\gg$.} \thmlab\semisimplecase

{\it Proof of Theorem \homospace~for semisimple $G$:} First we assume $G$ is semisimple. Since $G$ acts almost effectively on $G/H$, $\gh$ does not contain any simple component of $\gg$, so the conclusion of Lemma \surjectivity~holds. Fix a basis $\{\xi_1,\dots,\xi_n\}$ of $\gg$ and a corresponding dual basis $\{\xi'_1,\dots,\xi'_n\}$ of $\gg^*$ (relative to the Killing form), such that $\xi_1,\dots,\xi_h$ is a basis of $\gh$. By Lemma \vanishing, it suffices to construct a $G$-invariant, $G$-chiral horizontal element $\alpha\in \cW(\gg)\otimes \cQ'(G/H)$ satisfying $L^{tot}_{\xi}\circ_1\alpha = \beta^{\xi}\otimes 1$ for all $\xi\in\gg$.

In order to study $G/H$ as a $G$-space under left multiplication, it is convenient to regard $G$ as a $G\times H$-space, on which $G$ acts on the left and $H$ acts on the right. The right $H$-action induces compatible actions of $H$ and $\gs\gh[t]$ on $\cQ'(G)$ which commute with the actions of $G$ and $\gs\gg[t]$ coming from the left $G$-action. By Lemma 3.9 of \LLSI, the projection $\pi:G\ra G/H$ induces an isomorphism of vertex algebras $$\pi^*:\cQ'(G/H)) \ra \cQ'(G)_{H-bas}.$$ Moreover, by declaring that $H$ and $\gs\gh[t]$ act trivially on $\cW(\gg)$, we may extend the actions of $H$ and $\gs\gh[t]$ to $\cW(\gg)\otimes \cQ'(G)$. We identify the complexes $\cW(\gg)\otimes \cQ'(G/H)$ and $\cW(\gg)\otimes \cQ'(G)_{H-bas}$ and regard $\cW(\gg)\otimes \cQ'(G/H)_{H-bas}$ as a subcomplex of $\cW(\gg)\otimes \cQ'(G)$.
Thus in order to prove Theorem \homospace, it suffices to find a $G$-invariant, $G$-chiral horizontal element $\alpha\in\cW(\gg)\otimes \cQ'(G)_{H-bas}$ satisfying $L^{tot}_{\xi}\circ_1\alpha = \beta^{\xi}\otimes 1$ for all $\xi\in\gg$. In order to deal with all the operators $L^{tot}_{\xi}\circ_1$ simultaneously, it is convenient to define a new operator $$\cL:\cW(\gg)\otimes \cQ'(G)\ra \gg\otimes\cW(\gg)\otimes\cQ'(G)$$ sending $\omega\mapsto \xi_k\otimes L^{tot}_{\xi_k}\circ_1\omega$. Clearly $\cL$ is $G$-equivariant, and the condition $L^{tot}_{\xi}\circ_1\alpha = \beta^{\xi}\otimes 1$ for all $\xi\in\gg$ is equivalent to $\cL(\alpha) = \xi_k\otimes \beta^{\xi_k}\otimes 1$.

Let $f\in \cC^{\infty}(G)$ be a smooth function, and fix $\zeta\in\gg$, $\xi\in\gg$, and $\eta'\in\gg^*$. Then for $k=1,\dots,n$ we have
\eqn\circone{L^{tot}_{\xi_k}\circ_1 (\beta^{\zeta}b^{\xi} c^{\eta'}\otimes f) = L^{\cW}_{\xi_k}\circ_1 (\beta^{\zeta}b^{\xi} c^{\eta'}\otimes f) =  \beta^{\zeta}\otimes \bra [\xi_k,\xi],\eta'\ket f = \beta^{\zeta}\otimes \bra \xi_k, ad^*_{\xi}\eta'\ket f .} 

By Lemma \surjectivity, there exist elements $\chi_i\in \gg$ and $\eta'_i\in (\gg/\gh)^*$ for which $ad^*_{\chi_i}\eta'_i = \xi'_i$, for $i=1,\dots,n$. Then
$$L^{tot}_{\xi_k} \circ_1(\sum_i\beta^{\xi_i} b^{\chi_i} c^{\eta'_i}\otimes 1) =  \beta^{\xi_i}\otimes \bra \xi_k,\xi_i'\ket =  \beta^{\xi_k}\otimes 1,$$ so that
\eqn\condition{\cL(\sum_i\beta^{\xi_i} b^{\chi_i} c^{\eta'_i}\otimes 1) = \xi_k\otimes\beta^{\xi_k}\otimes 1.}
However, $\sum_i\beta^{\xi_i} b^{\chi_i} c^{\eta'_i}\otimes 1$ is not $G$-invariant. We seek a $G$-invariant element 
$$\alpha_0 = \sum_{j,k,l} \beta^{\xi_j} b^{\xi_k} c^{\xi'_l}\otimes f_{jkl}$$ which also satisfies \condition.

We will construct $\alpha_0$ using the connections one-forms coming from both the left and right actions of $G$ on itself, which we denote by $\theta^{\xi'}$, $\bar{\theta}^{\xi'}$, respectively, for $\xi'\in\gg^*$. We denote the $\gs\gg[t]$-algebra structure on $\cQ'(G)$ coming from the {\it right} $G$-action by $(\xi,\eta)\mapsto \bar{L}_{\xi}+\bar{\iota}_{\eta}$. Evaluating the functions $\iota_{\xi}\circ_0\bar{\theta}^{\xi'}$ and $\bar{\iota}_{\xi}\circ_0\theta^{\xi'}$ at the identity $e\in G$, we have
\eqn\evalid{\iota_{\xi}\circ_0 \bar{\theta}^{\xi'}|_e = \bra \xi,\xi'\ket=\bar{\iota}_{\xi}\circ_0 \theta^{\xi'} |_e.}
Define $$\alpha_0 = \sum_{i,j,k,l}\beta^{\xi_j} b^{\xi_k} c^{\xi'_l}\otimes \bar{\iota}_{\xi_i} (\theta^{\xi'_j})\bar{\iota}_{\chi_i}(\theta^{\xi'_k}) \iota_{\xi_l} (\bar{\theta}^{\eta'_i}).$$ Clearly $\alpha_0$ is $G$-invariant, and $ \alpha_0|_e = \sum_i\beta^{\xi_i} b^{\chi_i} c^{\eta'_i}\otimes 1$, by \evalid. Acting by $\cL$ we see that
$(\cL (\alpha_0))|_e = \xi_k\otimes \beta^{\xi_k}\otimes 1$. Finally, since $\alpha_0$ is $G$-invariant and the operator $\cL$ is $G$-equivariant, it follows that $\cL(\alpha_0) = \xi_k\otimes\beta^{\xi_k}\otimes 1$ at every point of $G$, as desired.

Our next step is to correct $\alpha_0$ to make it $G$-chiral horizontal without destroying $G$-invariance or condition \condition. Note that for $r \geq 0$,
$$\iota^{tot}_{\xi_t}\circ_r\alpha_0 = b^{\xi_t}\circ_r\alpha_0 = \delta_{r,0} \sum_{i,j,k}\beta^{\xi_j} b^{\xi_k} \otimes \bar{\iota}_{\xi_i} (\theta^{\xi'_j})\bar{\iota}_{\chi_i}(\theta^{\xi'_k}) \iota_{\xi_t} (\bar{\theta}^{\eta'_i}).$$
Let $$\alpha_1 = - \sum_{i,j,k,l}\beta^{\xi_j} b^{\xi_k}\otimes \bar{\iota}_{\xi_i}( \theta^{\xi'_j}) \bar{\iota}_{\chi_i}(\theta^{\xi'_k}) \iota_{\xi_l} (\bar{\theta}^{\eta'_i}) \theta^{\xi'_l}.$$ An OPE calculation shows that for $r\geq 0$
\eqn\alphaone{L^{tot}_{\xi_t}\circ_r\alpha_1 = 0,~~~~\iota^{tot}_{\xi_t} \circ_r \alpha_1 = -\delta_{r,0} \sum_{i,j,k}\beta^{\xi_j} b^{\xi_k} \otimes \bar{\iota}_{\xi_i} (\theta^{\xi'_j})\bar{\iota}_{\chi_i}(\theta^{\xi'_k}) \iota_{\xi_t} (\bar{\theta}^{\eta'_i}).} Let $\alpha = \alpha_0+\alpha_1$. It follows from \alphaone~that $\alpha$ is $G$-invariant, $G$-chiral horizontal, and satisfies $\cL(\alpha)=\xi_k\otimes\beta^{\xi_k}\otimes 1$.

We need to correct $\alpha$ so that it lies in $\cW(\gg)\otimes \cQ'(G)_{H-bas}$, without destroying the above properties. First, we claim that $\alpha_0$ is already $H$-chiral horizontal. This is clear since $\alpha$ is a sum of terms of the form $\beta^{\xi_j} b^{\xi_k} c^{\xi'_l}\otimes f_{jkl}$ where $f_{jkl}\in \cC^{\infty}(G)$, and $\bar{\iota}_{\xi}\circ_r$ lowers degree and only acts on the second factor of $\cW(\gg)\otimes \cQ'(G)$ for $\xi\in\gh$. 

Second, we claim that $\alpha_1$ is $H$-chiral horizontal as well. First note that for $\xi\in\gh$, $\bar{\iota}_{\xi}\circ_r$ acts by derivations on $\cQ'(G)$ for all $r\geq 0$, so it suffices to show that it acts by zero on each term of the form $\bar{\iota}_{\xi_i}( \theta^{\xi'_j})$, $ \bar{\iota}_{\chi_i}(\theta^{\xi'_k})$ and  $\iota_{\xi_l} (\bar{\theta}^{\eta'_i}) \theta^{\xi'_l}$. Clearly $\bar{\iota}_{\xi}\circ_r$ acts by zero on $\bar{\iota}_{\xi_i}( \theta^{\xi'_j})$ and $ \bar{\iota}_{\chi_i}(\theta^{\xi'_k})$ since these terms have degree 0 and $\bar{\iota}_{\xi}\circ_r$ has degree $-1$. Next, note that for each $\eta'_i\in (\gg/\gh)^*$ we have
$$\iota_{\xi_l}(\bar{\theta}^{\eta'_i})\theta^{\xi'_l} = \bar{\theta}^{\eta'_i}.$$ which can be checked by applying $\iota_{\eta}$, $\eta\in\gg/\gh$ to both sides. Since $\bar{\iota}_{\xi}\circ_r \bar{\theta}^{\eta'_i} = 0$ for all $\xi\in\gh$ and $\eta'_i\in (\gg/\gh)^*$, the claim follows.

Finally, if $\alpha$ is not $H$-invariant, we can make it $H$-invariant by averaging it over $H$, that is, we take $\alpha' = \frac{1}{|H|} \int_H h \alpha~d\mu$, where $d\mu$ is the Haar measure on $H$. Since the $G$ and $H$ actions commute, $\alpha'$ is $G$-invariant and $G$-chiral horizontal, and $\cL(\alpha') = \xi_k\otimes \beta^{\xi_k}\otimes 1$.  Moreover, since $\alpha$ is $H$-chiral horizontal, it follows that $h\alpha$ is still $H$-chiral horizontal for all $h\in H$, so $\alpha'$ is $H$-chiral horizontal as well. Since $\alpha'$ is $H$-invariant and lives in $\cW(\gg)[2]\otimes \cQ'(G)[0]$, $\alpha'$ is in fact $H$-chiral invariant, so that $\alpha'\in\cW(\gg)\otimes \cQ'(G)_{H-bas}$, as desired. $\Box$

{\it Proof of Theorem \homospace~for general $G$:} As in Lemma \finitecover, let $\tilde{G}$ be a finite cover of $G$ of the form $K\times T$, where $K$ is semisimple and $T$ is a torus. Then 
$$\H^*_G(\cQ'(G/H)) \cong \H^*_{\tilde{G}}(\cQ'(G/H)) \cong \H^*_{\tilde{G}}(\cQ'(\tilde{G}/\tilde{H})).$$ Here $\tilde{H}$ is the inverse image of $H$ under the projection $\tilde{G}\ra G$, which is a finite cover of $H$. Hence without loss of generality we may assume that $G$ is already of the form $K\times T$.

Since $G$ acts almost effectively on $G/H$, $H\cap T$ is finite, and since $T$ and $H$ commute, $H$ is a normal subgroup of $HT$. The group $HT/H$ is a torus, which we denote by $T'$. The natural map $T\ra T'$ is a finite cover, and we have a principal $T'$-bundle $G/H\ra G/HT$.

Since $K$ is semisimple and $K$ acts almost effectively on $$G/HT\cong (K\times T)/HT \cong K/(HT\cap K),$$ we can apply Theorem \homospace~in the semisimple case; there exists a $K$-chiral horizontal, $K$-invariant element $\alpha\in\cW(\gk)\otimes\cQ'(G/HT)$ such that $L_\xi^{tot}\circ_1\alpha=\beta^\xi\otimes 1$ for $\xi\in\gk$.

Since $G/H$ is a $T'$-bundle over $G/HT$, the projection $G/H\ra G/HT$ induces an isomorphism $\cQ'(G/HT)\ra \cQ'(G/H)_{T'-bas}$. Hence we have an injection $\cW(\gk)\otimes \cQ'(G/HT)\hookrightarrow \cW(\gg)\otimes \cQ'(G/H)$. Let $\alpha_K$ denote the image of $\alpha$ under this map. Clearly $\alpha_K$ is invariant under $L_{\eta}^{G/H}\circ_r$ and $\iota_{\eta}^{G/H}\circ_r$ for $\eta\in \gt$ and $r\geq 0$ because $\alpha_K$ is a sum of terms of the form $\alpha_i\otimes \omega_i$ with $\omega_i\in\cQ'(G/H)_{T'-bas}$. Similarly, $L^{\cW}_{\eta}\circ_r$ and $\iota_{\eta}^{\cW}\circ_r$ act trivially on $\alpha_K$ for $\eta\in\gt$ and $r\geq 0$ since each $\alpha_i$ does not depend on $b^{\eta},c^{\eta'},\beta^{\eta},\gamma^{\eta'}$ for $\eta\in\gt', \eta'\in(\gt')^*$. Hence $\alpha_K$ is $T'$-chiral basic. If $\alpha_K$ is not $T$-invariant, we can make it $T$-invariant by averaging it over the finite group $\Gamma = Ker (T\ra T')$.

Clearly $\alpha_K$ is $G$-chiral horizontal and $G$-invariant, and satisfies $L^{tot}_{\xi}\circ_1 \alpha_K = \beta^{\xi}\otimes 1$ for $\xi\in \gk$ and $L^{tot}_{\eta}\circ_1\alpha_K = 0$ for $\eta\in\gt$. We will construct another element $\alpha_T\in\cW(\gg)\otimes\cQ'(G/H)$ which is $G$-invariant, $G$-chiral horizontal, and satisfies $L^{tot}_{\xi}\circ_1 \alpha_T = 0$ for $\xi\in\gk$ and $L^{tot}_{\eta}\circ_1 \alpha_T = \beta^{\eta}\otimes 1$ for $\eta\in\gt$. Then $\alpha_K+\alpha_T$ satisfies all the conditions of Lemma \vanishing.

Since $H\cap T$ is finite, the action of $T$ on $G/H$ is locally free. Then there are connection forms $\theta^{\eta'}\in\Omega^1(G/H)$ satisfying $\iota_{\eta}\theta^{\eta'} = \bra\eta',\eta \ket$. Let $\Gamma^{\eta'} = g\circ_0 \theta^{\eta'}$, which has degree zero and weight one. Let $\eta_i$ be an orthonormal basis of $\gt$ (relative to the Killing form), and define
$$\alpha = \beta^{\eta_i}\otimes \Gamma^{\eta'_i}.$$ Clearly $\alpha$ is $T$-invariant, $T$-chiral horizontal, and satisfies $L^{tot}_{\eta}\circ_1 \alpha = \beta^{\eta}$ for $\eta\in \gt$. Since $K$ and $T$ commute, we can average $\alpha$ over $K$ without destroying the above properties to make it $G$-invariant. 

We claim that $\alpha$ is in fact $G$-chiral horizontal. It is enough to check this in local coordinates $\gamma^i$ on $G/H$. For $\xi'\in\gk^*$, $\Gamma^{\xi'}$ is of the form $f^k\partial\gamma^k$, where the $f^k$ are smooth functions. Similarly, for each $\xi\in\gk$, the vertex operator $\iota^{G/H}_{\xi}$ is locally of the form $f^kb^k$. Hence $\Gamma^{\xi'}$ commutes with $\iota^{G/H}_{\xi}$. Since $\alpha$ does not depend on $c^{\xi'}$ for $\xi'\in\gk^*$, it follows that $\alpha$ commutes with $b^{\xi}$ for $\xi\in\gk$. Hence $\iota^{tot}_{\xi} = b^{\xi} + \iota^{G/H}_{\xi}$ commutes with $\alpha$ for $\xi\in\gk$, as claimed. 

However, $\alpha$ need not satisfy $L^{tot}_{\xi}\circ_1\alpha = 0$ for $\xi\in \gk$. The last step is to correct $\alpha$ so that this property holds. Let $$\alpha_T = \alpha + \theta^{\xi_i}_{\cS} L_{\xi_i}^{tot}\circ_1 \alpha,$$ where $i$ runs over a basis of $\gk$. An OPE calculation using the fact that $L^{\cW}_{\xi_k}\circ_1 \theta^{\xi_i}_{\cS} = -\delta_{i,k}$ shows that $L^{tot}_{\xi}\circ_1 \alpha_T = 0$ for $\xi\in\gk$. Finally, since $L^{tot}_{\xi}\circ_1$ and $L^{tot}_{\eta}\circ_1$ commute for $\xi\in\gk$ and $\eta\in\gt$, it follows that for all $\eta\in\gt$ we have
$$L^{tot}_{\eta}\circ_1 (\theta^{\xi_i}_{\cS} L_{\xi_i}^{tot}\circ_1 \alpha) = \theta^{\xi_i}_{\cS} L^{tot}_{\xi_i}\circ_1 (L^{tot}_{\eta}\circ_1\alpha) =  \theta^{\xi_i}_{\cS} L^{tot}_{\xi_i}\circ_1 (\beta^{\eta}\otimes 1) = 0.
$$ Thus $L^{tot}_{\eta}\circ_1 \alpha_T = L^{tot}_{\eta}\circ_1 \alpha = \beta^{\eta}\otimes 1$, as desired. $\Box$

\subsec{Finite-dimensionality of $\H^*_G(\cQ'(M))$ for compact $M$}
The subspace $\cW(\gg)^p[n]\subset \cW(\gg)$ of degree $p$ and weight $n$ is finite-dimensional, so $\H^p_G(\C)[n]$ is finite-dimensional for all $p\in\Z$ and $n\geq 0$. Similarly, since $G/H$ is compact for any closed subgroup $H\subset G$, Theorem \homospace~implies that $\H^p_G(\cQ'(G/H))[n]$ is finite-dimensional for all $p,n$. In this section, we show that if we replace $G/H$ with an arbitrary compact $M$, $\H^p_G(\cQ'(M))[n]$ is finite-dimensional for all $p,n$ as well. This generalizes a well-known classical result in the case $n=0$. Hence the generating function 
$$\chi(G,M) = \sum_{p,n} dim~\H^p_G(\cQ'(M))[n]~z^p q^n$$ is a well-defined invariant of the $G$-manifold $M$. 
\lemma{If $M$ has a finite cover $\{U_1, \dots ,U_m\}$ of $G$-invariant open sets, such that
$$
dim~\H_G^p(\cQ'(U_{i_1}\cap\cdots\cap U_{i_k}))[n]<\infty,
$$
for each $p,n$ and for fixed indices $i_1,\dots ,i_k$, then $\H^p(\cQ'(M))[n]$ is finite-dimensional.}\thmlab\opencover
\proof
This is the standard generalized Mayer-Vietoris argument by induction on $m$. If $m=1$, there is nothing to prove. For $m=2$, this is the usual Mayer-Vietoris argument. Put $V=U_1\cup\cdots\cup U_{m-1}$. Then we have
$$\cdots\ra \H_G^{p-1}(\cQ'(V\cap U_m))\ra\H_G^p(\cQ'(M))\ra\H_G^p(\cQ'(V))\oplus\H_G^p(\cQ'(U_m))\ra\cdots,$$
which we restrict to a given weight $n$. By inductive hypothesis, $\H_G^p(\cQ'(V))[n]$ and $\H_G^p(\cQ'(U_m))[n]$ are finite-dimensional, so the third term is finite dimensional. Note that $V\cap U_m$ is covered by the open sets $U_i\cap U_m$, $i=1,\dots ,m-1$, and their multiple intersections also have finite-dimensional cohomology at fixed $p,n$. Since there are $m-1$ such open sets, the inductive hypothesis can be applied again. Thus the first term is also finite-dimensional, implying that the second term is finite-dimensional as well. $\Box$

\lemma{Suppose the $G$-manifold $M$ is a fiber bundle whose general fiber is $G/H$. If $M$ is compact, then $\H_G^p(\cQ'(M))[n]$ is finite-dimensional.}
\proof Choose a local trivializing cover of the bundle, which we can further refine to a good cover on the base $M/G$, i.e., each multiple intersection of open sets is a ball $B$. The preimage under $M\ra M/G$ of each multiple intersection of the open sets is equivariantly diffeomorphic to $G/H\times B$. So we can cover $M$ by finitely many open sets whose multiple intersections are equivariantly contractible to $G/H$. Since $\H^p_G(\cQ'(G/H))[n]$ is finite-dimensional for each $p,n$, the claim follows by the preceding lemma. $\Box$

Given a closed subgroup $H\subset G$, let $M_{(H)}$ denote the subset of $M$ consisting of points with isotropy group conjugate to $H$. $M_{(H)}$ is a closed submanifold of $M$, which may be regarded as a $G/H$-fiber bundle over the manifold $M_{(H)}/G$. By the preceding lemma, $\H^p_G(\cQ'(M_{(H)}))[n]$ is finite-dimensional for each $p,n$.
\theorem{ Suppose $M$ is compact. Then $\H^p_G(\cQ'(M))[n]$ is finite-dimensional.}\thmlab\finitedimension
\proof
Since $M$ is compact, there are only finitely many conjugacy classes $(H)$ for which $M_{(H)}$ is nonempty. For $dim~H>0$, each such $M_{(H)}$ has a $G$-invariant tubular neighborhood $U_{(H)}$ which is equivariantly contractible to $M_{(H)}$. $M$ has a finite cover consisting of the $U_{(H)}$ together with the open set $U$ of points with finite isotropy group. Without loss of generality we can shrink each $U_{(H)}$ so it contains only two orbit types: $(H)$ and $(e)$. 

By homotopy invariance, $\H^*_G(\cQ'(U_{(H)})) = \H^*_G(\cQ'(M_{(H)}))$, which is finite-dimensional at each $p,n$. The action of $G$ on $U$ is locally free so $\H^*_G(\cQ'(U))_+ = 0$. Since the multiple intersections of the $U_{(H)}$ all lie in $U$, they also have no higher-weight cohomology. It follows from Lemma \opencover~that for $n>0$, $\H^p_G(\cQ'(M))[n]$ is finite-dimensional for all $p$. For $n=0$, the finite-dimensionality of $\H^*_G(\cQ'(M))[0] = H^p_G(M)$ is classical. $\Box$

\newsec{The Structure of $\H^*_G(\cQ(M))$}

For any $G$-manifold $M$, Theorem \Qstructure~gives a complete description of $\H^*_G(\cQ(M))$ relative to the family of vertex algebras $\H^*_K(\C)$ for connected normal subgroups $K\subset G$. In particular, this result shows that $\H^*_G(\cQ(M))_+$ is only sensitive to $Ker(G\ra Diff(M))$, and in contrast to $\H^*_G(\cQ'(M))_+$, it carries no other geometric information about $M$. 

{\it Proof of Theorem \Qstructure:} We may assume that $G$ acts almost effectively on $M$. As usual, we need to show that the operator $L^{tot}\circ_1$ coming from the quasi-conformal structure on $\H^*_G(\cQ(M))$ acts by zero. By Lemma \vanishing, it is enough to construct a chiral horizontal, $G$-invariant element $\alpha\in \cW(\gg)\otimes \cQ(M)$ for which $L^{tot}_{\xi}\circ_1\alpha = \beta^{\xi}\otimes 1$ for all $\xi\in \gg$. Since $M$ can be covered by $G$-invariant open sets of the form $G\times_H V$ for some closed subgroup $H\subset G$ and some $H$-module $V$, it is enough to construct $\alpha\in \cW(\gg)\otimes \cQ(G\times_H V)$ for any $H$ and $V$, by Lemma \locI~and Remark \locII. Unlike the functor $\H^*_G(\cQ'(-))$, $\H^*_G(\cQ(-))$ is not a homotopy invariant, and in general $\H^*_G(\cQ(G\times_H V)) \neq \H^*_G(\cQ(G/H))$.

Fix $M = G\times_H V$, and assume that $G$ acts almost effectively on $M$. Note that the action of $G$ on the zero-section $G/H\subset M$ need not be almost effective. Let $K = Ker (G\ra Diff (G/H))$, and suppose first that $K$ is finite. Then by Theorem \homospace, there exists $\alpha\in\cW(\gg)\otimes \cQ'(G/H)$ which we can pull back to $\alpha\in\cW(\gg)\otimes \cQ'(M)\subset \cW(\gg)\otimes \cQ(M)$ with the same properties. 

So assume that $K$ has positive dimension, and let $K_0$ denote the identity component. By Lemma \finitecover, we may assume that $G$ is of the form $G_1\times\cdots\times G_k\times T$, where the $G_i$ are simple and $T$ is a torus. Any connected normal subgroup of $G$ splits, so $K_0$ is a product of a subset of the $G_i$ and a subtorus $T'\subset T$. Hence we may write $G = K_0\times N$ for some $N$. Likewise, since $K_0$ is also a normal subgroup of $H$, we may write $H = K_0\times L$. Then
$$M = (K_0\times N)\times_{K_0\times L} V = N\times_L V,$$ where $N$ acts on the left and $K_0$ acts only on the factor $V$, and the actions of $K_0$ and $L$ on $V$ commute. Let $\gn$ denote the Lie algebra of $N$. Clearly $Ker(N\ra Diff(N/L))$ is finite, so by the proof of Theorem \homospace~there exists an $N$-invariant, $N$-chiral horizontal element $\alpha_N\in\cW(\gn)\otimes \cQ'(N/L)$ such that $$L^{tot}_{\xi}\circ_1\alpha_N = \beta^{\xi}\otimes 1$$ for $\xi\in\gn$. Since $K_0$ acts trivially on $\cW(\gn)\otimes\cQ'(N/L)$, we can lift $\alpha_N$ to a $G$-chiral horizontal $G$-invariant element, also denoted by $\alpha_N$, in $\cW(\gn)\otimes\cQ'(N\times_L V)\subset\cW(\gg)\otimes\cQ'(N\times_LV)$ via the map induced by the projection map $N\times_L V\ra N/L$ such that $L^{tot}_\xi\circ_1\alpha_N=\beta^{\xi}\otimes 1$ for $\xi\in\gn$. We can view this element as lying in $\cW(\gg)\otimes\cQ(N\times_L V)$.

On the other hand, the linear action of $K_0$ on $V$ is faithful because the action of $G$ on $M$ is almost effective. Thus there exists a $K_0$-chiral horizontal, $K_0$-invariant element $\alpha_K\in\cW(\gk)\otimes\cQ(V)$ such that $L^{tot}_\xi\circ_1\alpha_K=\beta^\xi\otimes 1$ for $\xi\in\gk$. Now $L$ acts only on the $\cQ(V)$ factor and it commutes with $K_0$ action. Thus by averaging over $L$, we may assume that $\alpha_K$ is $L$-invariant, and still satisfies the same equation.

Recall from \linearspace~that $\alpha_K = \beta^{\xi_i}\otimes \Gamma^{\xi_i} - \beta^{\xi_i} b^{\xi_j}\otimes \iota_{\xi_j}\circ_0 \Gamma^{\xi_i}$, where $\Gamma^{\xi_i} = \beta^{\rho(\xi_i)(x_k)}\gamma^{x'_k}$. Clearly $\Gamma^{\xi_i}$ is itself $L$-invariant (because the action of $L$ on $x_i$ and dual action on $x_i'$ amounts to only a change of basis in $V$.) Now observe that any $L$-invariant element of $\cQ(V)$ can be regarded as a global section in $\cQ(N\times_L V)$. (This is clear by covering the bundle by open sets $U\times V$ where $U\subset N/L$, and the observing that the transition functions are elements of $L$.) The element $\alpha_K$ is now $N$-chiral horizontal and $N$-invariant in $\cW(\gg)\otimes\cQ(N\times_L V)$ and satisfies $L^{tot}_\xi\circ_1\alpha_K=0$ for $\xi\in\gn$.

Finally, $\alpha=\alpha_N+\alpha_K$ is $G$-chiral horizontal, $G$-invariant, and satisfies
$$
L^{tot}_\xi\circ_1\alpha=\beta^\xi\otimes 1
$$
for all $\xi\in\gg=\gk\oplus \gn$. This completes the proof that $\H^*_G(\cQ(M))_+=0$. $\Box$

\subsec{The ideal property of $\H^*_G(\cQ(M))_+$ and $\H^*_G(\cQ'(M))_+$}
An {\it ideal} in a vertex algebra $\cA$ is a linear subspace which is closed under $\alpha\circ_n$ and $\circ_n\alpha$ for all $n\in\Z$ and $\alpha\in\cA$. As usual, if $f:\cA\ra \cB$ is a vertex algebra homomorphism and $\cI\in\cB$ is an ideal, then $f^{-1}(\cI)$ is an ideal of $\cA$.

\lemma {For any $G$, $\H^*_G(\C)_+$ is an ideal of $\H^*_G(\C)$.}
\proof Let $V$ be a faithful representation of $G$, and consider the chiral Chern-Weil map $\kappa_G:\H^*_G(\C)\ra\H^*_G(\cQ(V))$. At weight zero, this is the identity map $S(\gg^*)^G\ra S(\gg^*)^G$ and it vanishes beyond weight zero because $\H^*_G(\cQ(V))_+ = 0$. Hence $\H^*_G(\C)_+ = Ker(\kappa_G)$, which is clearly an ideal. $\Box$

In view of Theorem \Qstructure~and the preceding lemma, it is immediate that $\H^*_G(\cQ(M))_+$ is an ideal of $\H^*_G(\cQ(M))$ for any $M$.

\theorem{For any $G$ and $M$, $\H^*_G(\cQ'(M))_+$ is an ideal of $\H^*_G(\cQ'(M))$.}\thmlab\idealproperty

\proof Under the map $\phi:\H^*_G(\cQ'(M))\ra \H^*_G(\cQ(M))$ induced by the inclusion $\cQ'(M)\hookrightarrow\cQ(M)$, we have $\H^*_G(\cQ'(M))_+ = \phi^{-1}(\H^*_G(\cQ'(M))_+)$, which is an ideal since $\H^*_G(\cQ(M))_+$ is an ideal of $\H^*_G(\cQ(M))$. $\Box$

\newsec{The Structure of $\H^*_G(\cQ'(M))$}
In contrast to $\H^*_G(\cQ(M))$, $\H^*_G(\cQ'(M))$ typically contains nontrivial geometric information about $M$ beyond weight zero. Our goal in this section is to give a {\it relative} description of $\H^*_G(\cQ'(M))$ in terms of the vertex algebras $\H^*_K(\C)$ for connected normal subgroups $K\subset G$, together with certain geometric data about $M$. We will focus on three special cases: $G$ simple, $G=G_1\times G_2$ where $G_1,G_2$ are simple, and $G$ abelian. We first describe $\H^*_G(\cQ'(M))$ as a linear space using Mayer-Vietoris sequences together with Theorem \homospace, and then describe the vertex algebra structure of $\H^*_G(\cQ'(M))$.

Let $\cA = \bigoplus_{n\geq 0} \cA[n]$ and $\cB = \bigoplus_{n\geq 0} \cB[n]$ be weight graded vertex algebras, and let $f:\cA\ra \cB$ be a weight-preserving vertex algebra homomorphism. Assume that $\cA_+ = \bigoplus_{n>0} \cA[n]$ is an ideal of $\cA$.

\lemma{Suppose that the restriction of $f$ to $\cA_+$ is injective. Then the vertex algebra structure of $\cA$ is uniquely determined by the ring structure of $(\cA[0],\circ_{-1})$, the homomorphism $f$, and the vertex algebra structure of $\cB$.}\thmlab\vastruct

\proof We need to describe the map 
\eqn\circn{\circ_n: \cA[i]\otimes \cA[j]\ra \cA[i+j-n-1]} for all $i,j\geq 0$ and $n\leq i+j-1$ in terms of the above data. First, suppose that $n = i+j-1$. Since $\cA_+$ is an ideal, it follows that \circn~is zero unless $i=j=0$ and $n=-1$, which is known by hypothesis. So we may assume that $n<i+j-1$.

Let $a\in \cA[i]$ and $b\in\cA[j]$, and suppose that either $i>0$ or $j>0$. Since $\cA_+$ is an ideal, $a\circ_n b\in \cA_+$. Since $f:\cA_+\ra \cB_+$ is injective, $f$ is an isomorphism of $\cA_+$ onto its image $\cB'_+\subset\cB_+$, and the inverse map $f^{-1}:\cB'_+\ra \cA_+$ is well-defined. It follows that 
$$a\circ_n b = f^{-1}(f(a)\circ_n f(b)),$$ which is determined by $f$ and the vertex algebra structure of $\cB$.

Finally, suppose that $i=j=0$ and $n\leq -2$. Then $a\circ_n b\in\cA_+$, so as above, we have $a\circ_n b = f^{-1}(f(a)\circ_n f(b))$. $\Box$

We will apply Lemma \vastruct~in the case where $\cA = \H^*_G(\cQ'(M))$ and $\cB$ is another vertex algebra whose structure is known, to determine the vertex algebra structure on $\H^*_G(\cQ'(M))$.

\subsec{The case where $G$ is simple}
Let $G$ be simple and let $M$ be a $G$-manifold. The inclusion $i:M^G\ra M$ induces a map \eqn\restmap{{\bf i}^*: \H^*_G(\cQ'(M))\ra \H^*_G(\cQ'(M^G))} on chiral equivariant cohomology, whose restriction to the weight zero subspace $\H^*_G(\cQ'(M))[0]$ coincides with the classical map $i^*:H^*_G(M)\ra H^*_G(M^G)$. 

{\it Proof of Theorem \simpleloc:}
We may assume that $M^G$ is non-empty. Let $U_0$ be a $G$-invariant tubular neighborhood of $M^G$ and let $U_1=M\setminus M^G$. It suffices to show that $\H^*_G(\cQ'(U_1))_+ = 0$ and $\H^*_G(\cQ'(U_0\cap U_1))_+ = 0$. In this case, we have $\H^*_G(\cQ'(M))_+ = \H^*_G(\cQ'(U_0\cup U_1))_+$, and since $\H^*_G(\cQ'(U_0))_+ \cong \H^*_G(\cQ'(M^G))_+$ by homotopy invariance, the claim follows from a Mayer-Vietoris argument.

For each point $x\in U_1$, the isotropy group $G_x$ has positive codimension in $G$ since $G$ is connected. Let $U_x$ be a $G$-invariant neighborhood of the orbit $Gx$, which we may take to be a vector bundle of the form $G\times_{G_x} V$ whose zero-section is $Gx$.

By Theorem \homospace, there exists a $G$-invariant, $G$-chiral horizontal element $\alpha_x\in \cW(\gg)\otimes \cQ'(G/G_x)$ satisfying \circleone. Via the projection $U_x\ra Gx$, this pulls back to an element $\alpha_{U_x}\in \cW(\gg)\otimes \cQ'(U_x)$ satisfying the same conditions. Using a $G$-invariant partition of unity as in Lemma \locI, we can glue the $\alpha_{U_x}$ together to obtain $\alpha\in \cW(\gg)\otimes \cQ'(U_1)$ satisfying these conditions as well. It follows that $\H^*_G(\cQ'(U_1))_+ = 0$. Finally, the same argument shows that $\H^*_G(\cQ'(U_0\cap U_1))_+ = 0$.

Next, let $\cA = \H^*_G(\cQ'(M))$, $\cB = \H^*_G(\cQ'(M^G))$, and let $f:\cA\ra \cB$ be the map ${\bf i}^*$ given by \restmap. The hypothesis of Lemma \vastruct~is satisfied, and the ring structure of $(\H^*_G(\cQ'(M))[0],\circ_{-1})$ coincides with the ring structure of $(H^*_G(M),\cup)$, which is classical. As a vertex algebra, $\H^*_G(\cQ'(M^G)) \cong \H^*_G(\C)\otimes H^*(M^G)$ where $H^*(M^G)$ is regarded as a vertex algebra in which all products except $\circ_{-1}$ are trivial. By Lemma \vastruct, this determines the vertex algebra structure of $\H^*_G(\cQ'(M))$ uniquely. 

Note that $\H^*_G(\cQ'(M))_+ = \H^*_G(\C)_+\otimes H^*(M^G)$ may alternatively be described as \eqn\tensorover{\H^*_G(\C)_+\otimes_{S(\gg^*)^G} H^*_G(M^G).} The tensor product above is defined by the obvious relation $(a\circ_{-1} b) \otimes c=a\otimes(b\circ_{-1}c)$. Given $\alpha\in \H^*_G(\C)_+$ and $\omega\in H^*_G(M)$, it follows from Theorem \simpleloc~that $$\kappa_G(\alpha)\circ_{-1} \omega = \alpha\otimes i^*(\omega)\in \H^*_G(\C)_+\otimes_{S(\gg^*)^G} H^*_G(M^G).$$ Similarly, given $\omega,\nu\in H^*(M^G)$, $\alpha\otimes\omega$ and $\alpha\otimes \nu$ lie in $\H^*_G(\C)\otimes H^*(M^G)$, and 
$$(\alpha\otimes \omega)\circ_{-1}(\alpha\otimes \nu) = (\alpha\circ_{-1} \alpha) \otimes (\omega\cup \nu).$$ Hence both the classical restriction map $i^*$ and the ring structure of $H^*(M^G)$ are encoded in the vertex algebra structure of $\H^*_G(\cQ'(M))$. $\Box$

\subsec{Chiral equivariant cohomology of spheres}

In this section we prove Theorem \mainresult. Let $G$ be a compact, connected Lie group. A $G$-manifold $M$ is said to be {\it equivariantly formal} if the spectral sequence of the fibration $$M\hookrightarrow (M\times E)/G\rightarrow E/G$$ collapses. If $M$ is equivariantly formal, the map $i^*: H^*_G(M)\ra H^*_G(M^G)$ is injective, and $H^*_G(M)\cong H^*(M)\otimes H^*_G(pt)$ as a module over $H^*_G(pt)$. Moreover, if $M$ is $G$-equivariantly formal, $M$ is also $K$-equivariantly for any closed subgroup $K\subset G$. The following result of \GKM~gives a useful criterion for equivariant formality. 

\theorem{(Goresky-Kottwitz-MacPherson) If the homology $H_*(M,\R)$ can be represented by $G$-invariant cycles, $M$ is equivariantly formal.}\thmlab\gkmtheorem

The next theorem we will need is an immediate consequence of results in \OI\OII~which describe the fixed-point sets of group actions on disks.

\theorem{(Oliver) Let $F$ be a finite $CW$-complex. If $G$ is semisimple, there exists a smooth action of $G$ on a closed disk $D$ with fixed point set $D^G$ having the homotopy type of $F$. If $G$ is a torus, there exists such an action if and only if $F$ is $\Z$-acyclic.}
\thmlab\oliver

If $G$ acts smoothly on an $n$-dimensional disk $D$, we may glue together two copies of $D$ along their boundaries to obtain a smooth action of $G$ on the sphere $S^n$. It is immediate from Theorem \gkmtheorem~that a $G$-sphere is equivariantly formal if and only if it has a $G$-fixed point. This allows us to construct compact $G$-manifolds which have the same classical equivariant cohomology, but distinct chiral equivariant cohomology.

{\it Proof of Theorem \mainresult:} Let $G$ be simple, and let $F$ be a $CW$-complex with $3$ zero-cells and no higher-dimensional cells. Choose an $n$-dimensional disk $D$ with a smooth $G$-action such that $D^G$ has the homotopy type of $F$. Let $S_0$ be the copy of $S^{n}$ obtained by gluing together two copies of $D$ along their boundaries. Note that each connected component $\cC$ of $D^G$ gives rise to one component of $S_0^G$ (if $\cC\cap \partial D \neq \varnothing$), or two components of $S_0^G$ (if $\cC\cap \partial D_0= \varnothing$).  Hence $3\leq c_0 \leq 6$, where $c_0$ is the number of components of $S_0^G$. 

Given $x\in S_0^G$, we can find a $G$-invariant ball $B_0\subset S_0$ containing $x$, which intersects exactly one component of $S_0^G$. By removing $B_0$ from two copies of $S_0$ and then gluing them together along their boundaries, we obtain an $n$-dimensional $G$-sphere $S_1$ such that $S_1^G$ has either $2c_0-2$ components (if $\partial B\cap S_0^G = \varnothing$) or $2c_0-1$ components (if $\partial B \cap  S_0^G \neq \varnothing$). We continue this process as follows. Assume that $n$-dimensional $G$-spheres $S_0, S_1,\dots,S_{i-1}$ have been defined. Let $B_{i-1}\subset S_{i-1}$ be a $G$-invariant ball intersecting exactly one component of $S_{i-1}^G$. Define $D_i = S_{i-1}\setminus B_{i-1}$, and let $S_i$ be the sphere obtained by gluing two copies of $D_i$ along their boundaries. We thus obtain a sequence of $n$-dimensional spheres $S_0,S_1,S_2,\cdots$ with smooth $G$-actions, such that the number of fixed-point components $c_0,c_1,c_2,\dots$ are all distinct. Since $c_i = dim~H^0(S_i^G)$, it is immediate from Theorem \simpleloc~that the vertex algebras $\H^*_G(\cQ'(S_i))$ are all distinct.

To complete the proof of Theorem \mainresult, it remains to show that the classical equivariant cohomology rings $H^*_G(S_i)$ are all isomorphic to $H^*_G(pt)[\omega]/(\omega^2)$. Let $T$ be a maximal torus of $G$, and let $W$ be the Weyl group, so that $H^*_G(S_i) = H^*_T(S_i)^W$.

By Theorem \oliver, $D_i^T$ is acyclic, and by a Mayer-Vietoris argument, $S_i^T$ has the homology type of a $k$-dimensional sphere for some $k<n$. In particular, $dim~H^*(S_i^T) = 2$. Next, we claim that $S_i^T$ is connected, so that $k\geq 1$. To see this, note first that $D_i^T$ must be connected, since $D_i^T$ is acyclic and contains more than two points. It suffices to prove that $D_i^T\cap \partial D_i \neq \varnothing$. Suppose that $D_i^T$ lies in the interior of $D_i$, and note that $S_i \setminus D_i$ must also contain exactly one $T$-fixed
point component. Given a point $x\in D_i^T$, choose a small $G$-invariant ball $B\subset D_i$
containing $x$. Since $x$ is not an isolated $T$-fixed point, we may assume (by
taking $B$ sufficiently small) that $D_i^T\cap \partial B \neq \varnothing$. But $S_i\setminus B$ is then a $G$-invariant disk with more than one $T$-fixed point component, which is impossible by Theorem 7.2.

Since $S_i$ is $G$-equivariantly formal (and hence $T$-equivariantly formal), the map \eqn\eqform{i^*:H^*_T(S_i)\ra H^*_T(S_i^T)= H^*(S_i^T)\otimes H^*_T(pt)} is injective, and $H^*_T(S_i)\cong H^*(S_i)\otimes H^*_T(pt)$. As an $H^*_T(pt)$-module, $H^*_T(S_i)$ has generators $1$ and $\omega$ of degrees $0$ and $n$, respectively. Since $dim~H^*(S_i^T) = 2$, \eqform~becomes an isomorphism after localization along $H^*_T(pt)$. Letting $k$ be the dimension of $S_i^T$, since $dim~H^k(S_i^T) = 1$, we must have $$i^*(\omega)\in H^k(S_i^T)\otimes H^*_T(pt).$$ It follows that $i^*(\omega^2) = 0$, so $\omega^2 = 0$, and we have a direct sum decomposition $$H^*_T(S_i) = H^*_T(pt)\oplus H^*_T(pt)\omega.$$ Next, the Weyl group $W$ must preserve the subspace $H^*_T(pt) \omega\subset H^*_T(S_i)$, so we have
$$H^*_G(S_i) = H^*_T(S_i)^W = H^*_T(pt)^W \oplus (H^*_T(pt)\omega)^W,$$ which is a free module over $H^*_T(pt)^W = H^*_G(pt)$ having generators of degrees $0$ and $n$. Thus $H^*_G(pt)\omega$ contains an element of degree $n$, which must be $\omega$, so $H^*_G(S_i) \cong H^*_G(pt)[\omega]/(\omega^2)$ for each $i\geq 1$, as claimed. $\Box$

Next, we give examples of {\it morphisms} $f:M\ra N$ in the category of compact $G$-manifolds (that is, smooth $G$-equivariant maps) for which $f^*:H^*_G(N)\ra H^*_G(M)$ is an isomorphism (over $\Z$), but $\H^*_G(\cQ'(M))\neq \H^*_G(\cQ'(N))$.

Let $D_0$ be an $n$-dimensional disk with a smooth $G$-action for which $\varnothing \neq D^G\neq D$, and let $S_0$ be the corresponding sphere, as above. Fix a point $p\in S_0^G$, and let $U$ be a $G$-invariant neighborhood of $p$ equipped with a smooth $G$-equivariant map $\phi:U\ra \R^n$, where $G$ acts linearly on $\R^n$, $\phi(p)= 0$, and $\phi$ maps the closure $\bar{U}$ diffeomorphically onto the disk $D_1 = \{x\in \R^n: |x|\leq 1\}$. We identify $D_1/\partial D_1$ with another copy of $S^{n}$, which we denote by $S_1$. Note that $S_1$ is a smooth $G$-manifold and $S_1^G$ is a sphere $S^k$ for some $0\leq k< n$. Let $q\in S_1^G$ be the point which corresponds to $\partial D_1$ under the projection $\pi: D_1\ra S_1$, and let $g = \pi\circ \phi:U\ra S_1$, which is smooth and $G$-equivariant. Note that $g$ does {\it not} extend to a smooth map $S_0\ra S_1$. However, we can construct a new function $f$ which agrees with $g$ in a neighborhood of $p\in U$, which extends smoothly to all of $S_0$.

\lemma{There exists a smooth, $G$-equivariant map $f:S_0\ra S_1$ such that $f=g$ in a neighborhood of $p\in U$, and $f(S_0 \setminus U) = q$. Moreover, $f$ is smoothly (but not equivariantly) homotopic to the identity map $S^{n}\ra S^{n}$, so $f$ induces isomorphisms in both singular and equivariant cohomology with $\Z$-coefficients.}

\proof As above, we identify $U$ with the open disk $D^{\circ}_1 = \{(x: |x| <1\}$. Let $U'\subset U$ be a $G$-invariant neighborhood of $\partial U$ of the form $\{x\in U:~1-\epsilon<|x|\leq 1\}$. Clearly $g(U')$ is a $G$-invariant neighborhood of $q$, which we identify with another open disk $D_2^{\circ} =  \{y:  |y| <1\}$ equipped with a linear action of $G$. Note that the radius $|y|$ is a $G$-invariant function on $D_2$.

Choose a smooth function $h:[0,1]\ra [0,1]$ such that $h(t) = 0$ for $0\leq t < \frac{1}{3}$, and $h(t) = 1$ for $2/3 < t \leq 1$. Define $f:U\ra S_1$ by $f(x) = g(x) h(|g(x)|)$ for $x\in U'$, and $f(x) = g(x)$ for $x\in U\setminus U'$. Since $h(|y|)$ is $G$-invariant, $f$ is $G$-equivariant and since $f = g$ on a neighborhood of $U\setminus U'$, $f$ is smooth. Moreover, since $f$ maps a neighborhood of $\partial U$ to $q$ (which we have identified with $0\in D_2$), $f$ extends to a smooth, $G$-equivariant map $S_0\ra S_1$ sending $S_0\setminus U\ra q$, as desired. Finally, the fact that $f$ is homotopic to $id:S^{n}\ra S^{n}$ is clear because $S_0\setminus U$ is smoothly contractible to a point in $S_0$. $\Box$

\theorem{Suppose that $F$ is a $CW$-complex consisting of $m$ zero cells and no higher-dimensional cells. Let $D$ be a disk with a smooth $G$-action such that $D^G$ has the homotopy type of $F$, and let $S_0$ and $S_1$ be as above. If $m\geq 3$, $\H^*_G(\cQ'(S_0))\neq \H^*_G(\cQ'(S_1)$.}

\proof Since $S_1^G$ is a $k$-dimensional sphere for $0\leq k <n$,  we have $1\leq dim~H^0(S_1^G) \leq 2$. Since $m\leq dim~H^0(S_0^G) \leq 2m$, we have $dim~H^0(S_0^G) > dim~H^0(S_1^G)$ for $m\geq 3$. The claim then follows from Theorem \simpleloc. $\Box$

Even though $f:S_0\ra S_1$ is homotopic to $id:S^{n}\ra S^{n}$, this result does not contradict Theorem \hominv~because there is no {\it equivariant} homotopy between $f$ and $id$. Thus unlike the classical equivariant cohomology, the functor $\H^*_G(\cQ'(-))$ can distinguish $G$-manifolds $M$ and $N$ which admit a $G$-equivariant map which is a homotopy equivalence, as long as $M$ and $N$ are not equivariantly homotopic.

\subsec{The case $G=G_1\times G_2$, where $G_1,G_2$ are simple}
For simple $G$, Theorem \simpleloc~describes $\H^*_G(\cQ'(M))$ in terms of $\H^*_G(\C)$ together with classical geometric data. In this section, we give a similar description of $\H^*_G(\cQ'(M))_+$ in the case $G = G_1\times G_2$, where $G_1$ and $G_2$ are simple groups. As in the case of simple $G$, we first describe $\H^*_G(\cQ'(M))_+$ as a linear space, and then describe the vertex algebra structure. 

\theorem {Let $G_1,G_2$ be simple and let $G=G_1\times G_2$. For any $G$-manifold $M$, $\H^*_G(\cQ'(M))_+$ is linearly isomorphic to
$$\H^*_{G_1}(\C)_+\otimes H^*_{G_2}(M^{G_1}) \bigoplus 
\H^*_{G_2}(\C)_+\otimes H^*_{G_1}(M^{G_2}) \bigoplus \H^*_{G_1}(\C)_+\otimes\H^*_{G_2}(\C)_+\otimes H^*(M^G).$$} \thmlab\productsimple

Let $U_1,U_2$ be $G$-invariant tubular neighborhoods of $M^{G_1}, M^{G_2}$, respectively.
If $x\notin U_1 \cup U_2$, its stabilizer $G_x$ contains neither $G_1$ nor $G_2$, so $\H^*_G(\cQ'(G/G_x))_+ = 0$. A Mayer-Vietoris argument then shows that $\H^*_G(\cQ'(M))_+ = \H^*_G(\cQ'(U_1\cup U_2))_+$. Since $U_1$, $U_2$, and $U_1\cap U_2$ are equivariantly contractible to $M^{G_1}$, $M^{G_2}$, and $M^G$, respectively, we can replace $\H^*_G(\cQ'(U_1))_+$, $\H^*_G(\cQ'(U_2))_+$, and $\H^*_G(\cQ'(U_1\cap U_2))_+$ with $\H^*_G(\cQ'(M^{G_1}))_+$, $\H^*_G(\cQ'(M^{G_2}))_+$, and $\H^*_G(\cQ'(M^{G}))_+$ in the Mayer-Vietoris sequence 
$$\cdots\ra\H^*_G(\cQ'(U_1\cup U_2))_+\ra \H^*_G(\cQ'(U_1))_+ \oplus\H^*_G(\cQ'(U_2))_+\ra \H^*_G(\cQ'(U_2\cap U_2))_+\ra\cdots,$$ obtaining
\eqn\dumb{\cdots\ra\H^*_G(\cQ'(M))_+\ra \H^*_G(\cQ'(M^{G_1}))_+ \oplus\H^*_G(\cQ'(M^{G_2}))_+\ra \H^*_G(\cQ'(M^G))_+\ra\cdots.}

\lemma{The map $\phi: \H^*_G(\cQ'(M^{G_1}))_+ \oplus\H^*_G(\cQ'(M^{G_2}))_+\ra \H^*_G(\cQ'(M^G))_+$ appearing in \dumb~is surjective. Hence $\H^*_G(\cQ'(M))_+$ is canonically isomorphic to $Ker(\phi)$.}
\proof
Let $i_1:M^{G}\ra M^{G_1}$, $i_2:M^{G}\ra M^{G_2}$, $j_1:M^{G_1}\ra M$, $j_2:M^{G_2}\ra M$ denote the obvious inclusion maps. First we need to describe each of the spaces $\H^*_G(\cQ'(M^{G_1}))_+$, $\H^*_G(\cQ'(M^{G_2}))_+$, and $\H^*_G(\cQ'(M^G))_+$. Since $G$ acts trivially on $M^G$, $\H^*_G(\cQ'(M^G))_+ = \H^*_G(\C)_+ \otimes H^*(M^G)$, which is isomorphic to $$ \H^*_{G_1}(\C)_+\otimes H^*_{G_2}(pt)\otimes H^*(M^G) \bigoplus  H^*_{G_1}(pt) \otimes \H^*_{G_2}(\C)_+ \otimes H^*(M^G) \bigoplus  \H^*_{G_1}(\C)_+ \otimes \H^*_{G_2}(\C)_+ \otimes H^*(M^G).$$

Similarly, since $G_1$ acts trivially on $M^{G_1}$, $\H^*_G(\cQ'(M^{G_1})) = \H^*_{G_1}(\C) \otimes \H^*_{G_2}(\cQ'(M^{G_1}))$. Hence $\H^*_G(\cQ'(M^{G_1}))_+$ is isomorphic to
$$\H^*_{G_1}(\C)_+ \otimes H^*_{G_2}(M^{G_1}) \bigoplus 
H^*_{G_1}(pt) \otimes \H^*_{G_2}(\cQ'(M^{G_1}))_+ \bigoplus \H^*_{G_1}(\C)_+ \otimes \H^*_{G_2}(\cQ'(M^{G_1}))_+.$$ Since $G_2$ is simple, $\H^*_{G_2}(\cQ'(M^{G_1}))_+ = \H^*_{G_2}(\C)_+ \otimes H^*(M^G)$ by Theorem \simpleloc. Hence 
$\H^*_{G}(\cQ'(M^{G_1}))_+$ is isomorphic to
$$ \H^*_{G_1}(\C)_+ \otimes H^*_{G_2}(M^{G_1}) \bigoplus 
H^*_{G_1}(pt) \otimes \H^*_{G_2}(\C)_+ \otimes H^*(M^G) \bigoplus \H^*_{G_1}(\C)_+ \otimes \H^*_{G_2}(\C)_+ \otimes H^*(M^G).$$ 

Interchanging the roles of $G_1$ and $G_2$, we see that $ \H^*_{G}(\cQ'(M^{G_2}))_+$ is isomorphic to
$$\H^*_{G_2}(\C)_+ \otimes H^*_{G_1}(M^{G_2}) \bigoplus 
H^*_{G_2}(pt) \otimes \H^*_{G_1}(\C)_+ \otimes H^*(M^G) \bigoplus \H^*_{G_2}(\C)_+ \otimes \H^*_{G_1}(\C)_+ \otimes H^*(M^G).$$

Next, we need to describe the restriction of $\phi$ to the various summands of $\H^*_G(\cQ'(M^{G_1}))_+\oplus \H^*_G(\cQ'(M^{G_2}))_+$. Clearly $\phi$ maps the summand $\H^*_{G_1}(\C)_+ \otimes H^*_{G_2}(M^{G_1})\subset \H^*_G(\cQ'(M^{G_1}))_+$ to $$\H^*_{G_1}(\C)_+\otimes H^*_{G_2}(pt)\otimes H^*(M^G)\subset \H^*_G(\cQ'(M^G))_+,$$
acting by $id\otimes i^*_1$. (For this to make sense, we need to identify $\H^*_{G_1}(\C)_+\otimes H^*_{G_2}(pt)\otimes H^*(M^G)$ with $\H^*_{G_1}(\C)_+\otimes H^*_{G_2}(pt)\otimes_{H^*_{G_2}(pt)} H_{G_2}^*(M^G)$, as in \tensorover). Also, $\phi$ acts by $id$ on the remaining summands of $\H^*_G(\cQ'(M^{G_1}))_+$.

Similarly, $\phi$ maps the summand  $\H^*_{G_2}(\C)_+ \otimes H^*_{G_1}(M^{G_2}) \subset \H^*_G(\cQ'(M^{G_2}))_+$ to $\H^*_{G_2}(\C)_+ \otimes H^*_{G_1}(pt)\otimes H^*(M^G)$, acting by $-id\otimes i^*_2$. Finally, $\phi$ acts by $-id$ the identity on the remaining summands of $\H^*_G(\cQ'(M^{G_2}))_+$. The surjectivity of $\phi$ is now apparent. $\Box$

{\it Proof of Theorem \productsimple~:} The following notation will be convenient. Since $\H^*_G(\cQ'(M^{G_1}))_+\oplus \H^*_G(\cQ'(M^{G_2}))_+$ decomposes as the direct sum of six subspaces
$$ \H^*_{G_1}(\C)_+ \otimes H^*_{G_2}(M^{G_1}) \bigoplus 
H^*_{G_1}(pt) \otimes \H^*_{G_2}(\C)_+ \otimes H^*(M^G) \bigoplus \H^*_{G_1}(\C)_+ \otimes \H^*_{G_2}(\C)_+ \otimes H^*(M^G)$$  $$\bigoplus \H^*_{G_2}(\C)_+ \otimes H^*_{G_1}(M^{G_2}) \bigoplus 
H^*_{G_2}(pt) \otimes \H^*_{G_1}(\C)_+ \otimes H^*(M^G) \bigoplus \H^*_{G_2}(\C)_+ \otimes \H^*_{G_1}(\C)_+ \otimes H^*(M^G),$$
an element $\omega\in \H^*_G(\cQ'(M^{G_1}))_+\oplus \H^*_G(\cQ'(M^{G_2}))_+$ can be written uniquely as a $6$-tuple $(\omega_1,\dots,\omega_6)$.

Let $\alpha =  \sum_i \alpha_i\otimes \omega_i$ be an arbitrary element of $\H^*_{G_1}(\C)_+\otimes H^*_{G_2}(M^{G_1}).$ Via $i^*_1:H^*_{G_2}(M^{G_1})\ra H^*_{G_2}(M^G)$, $\omega_i\mapsto \sum_j p_{ij}\otimes \nu_{ij}\in H^*_{G_2}(pt) \otimes H^*(M^G) = H^*_{G_2}(M^G)$. Let $$\tilde{\alpha} = \sum_{i,j} \alpha_i\otimes p_{ij}\otimes \nu_{ij}\in \H^*_{G_1}(\C)_+\otimes H^*_{G_2}(pt)\otimes H^*(M^G),$$ and note that $(\alpha,0,0,0,\tilde{\alpha},0)$ lies in $Ker(\phi)$ by construction. Since $\alpha$ was arbitrary, the assignment $\alpha\mapsto (\alpha,0,0,0,\tilde{\alpha},0)$  identifies $\H^*_{G_1}(\C)_+ \otimes H^*_{G_2}(M^{G_1})$ with a linear subspace of $Ker(\phi)$.

Interchanging the roles of $G_1$ and $G_2$, for any $\alpha\in \H^*_{G_2}(\C)_+\otimes H^*_{G_1}(M^{G_2})$ we can find $\tilde{\alpha}\in $ such that $(0,\tilde{\alpha},0,\alpha,0,0)\in Ker(\phi)$. Hence the assignment $\alpha\mapsto (0,\tilde{\alpha},0,\alpha,0,0)$ identifies $\H^*_{G_2}(\C)_+\otimes H^*_{G_1}(M^{G_2})$ with another subspace of $Ker(\phi)$. 

Finally, given $\alpha \in \H^*_{G_2}(\C)_+ \otimes \H^*_{G_1}(\C)_+ \otimes H^*(M^G)$, $(0,0,\alpha,0,0,\alpha)\in Ker(\phi)$, so the map $\alpha\mapsto (0,0,\alpha,0,0,\alpha)$ identifies $\H^*_{G_2}(\C)_+ \otimes \H^*_{G_1}(\C)_+ \otimes H^*(M^G)$ with another subspace of $Ker(\phi)$. Clearly these three subspaces of $Ker(\phi)$ intersect trivially, and account for all of $Ker(\phi)$. $\Box$

Next, we use Lemma \vastruct~to describe the vertex algebra structure of $\H^*_G(\cQ'(M))$. Taking $\cA = \H^*_G(\cQ'(M))$ and $\cB = \H^*_G(\cQ'(M^{G_1}))\oplus \H^*_G(\cQ'(M^{G_2}))$, and $f:\cA\ra \cB$ the map 
$$\H^*_G(\cQ'(M))\ra \H^*_G(\cQ'(M^{G_1}))\oplus \H^*_G(\cQ'(M^{G_2}))$$ appearing in \dumb~it is clear that the hypothesis of Lemma \vastruct~holds. The ring structure of $\H^*_G(\cQ'(M))[0] = H^*_G(M)$ is classical, and the vertex algebra structure of $\H^*_G(\cQ'(M^{G_1}))\oplus \H^*_G(\cQ'(M^{G_2}))$ may be described completely in terms of $\H^*_{G_1}(\C)$, $\H^*_{G_2}(\C)$ and classical data because of the identity
$$\H^*_G(\cQ'(M^{G_1}))\oplus \H^*_G(\cQ'(M^{G_2})) = \H^*_{G_1}(\C)\otimes \H^*_{G_2}(\cQ'(M^{G_1})) \oplus \H^*_{G_2}(\C)\otimes \H^*_{G_1}(\cQ'(M^{G_2})).$$ Since $G_1,G_2$ are simple, the vertex algebra structures of both $\H^*_{G_1}(\cQ'(M^{G_2}))$ and $ \H^*_{G_2}(\cQ'(M^{G_1}))$ are given by Theorem \simpleloc. By Lemma \vastruct, this uniquely determines the vertex algebra structure of $\H^*_G(\cQ'(M))$.

Finally, via the identification of $\H^*_G(\cQ'(M))_+$ with $$\H^*_{G_1}(\C)_+\otimes H^*_{G_2}(M^{G_1}) \bigoplus 
\H^*_{G_2}(\C)_+\otimes H^*_{G_1}(M^{G_2}) \bigoplus \H^*_{G_1}(\C)_+\otimes\H^*_{G_2}(\C)_+\otimes H^*(M^G)$$ given by Theorem \productsimple, we can now describe all circle products in $\H^*_G(\cQ'(M))= H^*_G(M) \oplus \H^*_G(\cQ'(M))_+$. For example,

\item{$\bullet$} Given $\alpha\otimes \omega\in \H^*_{G_1}(\C)_+\otimes H^*_{G_2}(M^{G_1})$ and $\eta\otimes \nu \in\H^*_{G_2}(\C)_+\otimes H^*_{G_1}(M^{G_2})$, $$(\alpha\otimes\omega)\circ_{-1}(\eta\otimes\nu) = (\alpha\circ_{-1}\eta) \otimes i_1^*(\alpha)\cup i_2^*(\nu)\in \H^*_{G_1}(\C)_+\otimes\H^*_{G_2}(\C)_+\otimes H^*(M^G).$$ 

\item{$\bullet$} Given $a \in H^*_G(M)$ and $\alpha\otimes\omega \in \H^*_{G_1}(\C)_+\otimes H^*_{G_2}(M^{G_1})$, $$a \circ_{-1}(\alpha \otimes \omega) = \alpha \otimes j_1^*(a)\cup \omega \in \H^*_{G_1}(\C)_+\otimes H^*_{G_2}(M^{G_1}).$$ 

\item{$\bullet$} Given $a \in H^*_G(M)$ and $\alpha\otimes\omega \in \H^*_{G_2}(\C)_+\otimes H^*_{G_1}(M^{G_2})$, $$a\circ_{-1}(\alpha \otimes \omega) = \alpha \otimes j_2^*(a)\cup \omega \in \H^*_{G_2}(\C)_+\otimes H^*_{G_1}(M^{G_2}).$$

Note that the maps $i^*_1$, $i_2^*$, $j^*_1$ and $j^*_2$, as well as the ring structures of $H^*_{G_1}(M^{G_2})$, $H^*_{G_2}(M^{G_1})$, and $H^*_G(M^G)$ are encoded in the vertex algebra structure of $\H^*_G(\cQ'(M))$. For general $G$, we expect $\H^*_G(\cQ'(M))$ to depend on the family of vertex algebras $\H^*_K(\C)$ for connected normal subgroups $K\subset G$ for which $M^K$ is nonempty, together with the rings $H^*_{G/K}(M^K)$ and all maps $H^*_{G/K} (M^K)\ra H^*_{G/K'} (M^{K'})$ when $K\subset K'$.

\subsec{The case where $G$ is a torus $T$}
In this section, we study $\H^*_G(\cQ'(M))$ in the case where $G$ is a torus $T$. Recall from \LLSI~that $\H^*_T(\cQ'(M))$ can be computed using the small chiral Weil complex $\cC = \bra \gamma,c\ket \otimes \cQ'(M)$, with differential $K(0)\otimes 1 + 1\otimes d_{\cQ}$. Here $\bra \gamma,c\ket $ is the subalgebra of $\cW(\gt)$ generated by the $\gamma^{\xi'_i},c^{\xi'_i}$, $\xi'\in\gt^*$, and $K(0)$ is the chiral Koszul differential. $\H^*_T(\cQ'(M))$ is an abelian vertex algebra, i.e., a supercommutative algebra equipped with a derivation $\partial$ of degree 0 and weight 1. 

First we consider the case $T=S^1$. As usual, the inclusion $i:M^{S^1}\ra M$ induces a vertex algebra homomorphism ${\bf i}^*: \H^*_{S^1}(\cQ'(M))\ra \H^*_{S^1}(\cQ'(M^{S^1}))$ whose restriction to weight zero coincides with the classical map $i^*:H^*_{S^1}(M)\ra H^*_{S^1}(M^{S^1})$. The next result is analogous to Theorem \simpleloc~for simple group actions.
 
\theorem{(Positive-weight localization for circle actions) For any $S^1$-manifold $M$, 
${\bf i}^*: \H^*_{S^1}(\cQ'(M))_+\ra \H^*_{S^1}(\cQ'(M^{S^1}))_+$ is an isomorphism of vertex algebra ideals. Hence \eqn\circloc{\H^*_{S^1}(\cQ'(M))_+ \cong \H^*_{S^1}(\C)_+\otimes H^*(M^{S^1}).}Moreover, both the ring structure of $H^*(M^{S^1})$ and the map $i^*:H^*_{S^1}(M)\ra H^*_{S^1}(M^{S^1})$ are encoded in the vertex algebra structure of $\H^*_{S^1}(\cQ'(M))$.} \thmlab\circleloc
\proof Every point in $M$ is either an $S^1$-fixed point or has a finite isotropy group. If $M^{S^1}$ is nonempty, fix an $S^1$-invariant tubular neighborhood $U$ of $M^{S^1}$, and let $V = M\setminus M^{S^1}$. Since $S^1$ acts locally freely on $V$, $\H^*_{S^1}(\cQ'(V))_+ = 0 = \H^*_{S^1}(\cQ'(U\cap V))_+$. By a Mayer-Vietoris argument, $\H^*_{S^1}(\cQ'(M))_+ \cong \H^*_{S^1}(\cQ'(U))_+$, and by homotopy invariance $\H^*_{S^1}(\cQ'(U))_+ \cong \H^*_{S^1}(\cQ'(M^{S^1}))_+$. Finally, since $S^1$ acts trivially on $M^{S^1}$, we have $\H^*_{S^1}(\cQ'(M^{S^1}))_+ \cong \H^*_{S^1}(\C)_+\otimes H^*(M^{S^1})$.  

As for the vertex algebra structure, taking $\cA = \H^*_{S^1}(\cQ'(M))$, $\cB = \H^*_{S^1}(\cQ'(M^{S^1}))$ and $f = {\bf i}^*$, the hypothesis of Lemma \vastruct~is clearly satisfied. The ring structure of $\H^*_{S^1}(\cQ'(M))[0] = H^*_{S^1}(M)$ is classical, and $\H^*_{S^1}(\cQ'(M^{S^1})) = \H^*_{S^1}(\C)\otimes H^*(M^{S^1}))$ as a vertex algebra. This determines the vertex algebra structure of $\H^*_{S^1}(\cQ'(M))$. As in the case where $G$ is simple, both the classical restriction map $i^*$ and the ring structure of $H^*(M^{S^1})$ are encoded in the vertex algebra structure of $\H^*_{S^1}(\cQ'(M))$. $\Box$

Recall from Theorem 6.1 of \LLI~that $\H^*_{S^1}(\C)$ is just the polynomial algebra $\C[\gamma,\partial\gamma,\partial^2\gamma,\dots]$. Hence $\H^*_{S^1}(\C)_+$ is the ideal $\bra \partial\gamma,\partial^2\gamma,\cdots\ket\subset \C[\gamma,\partial\gamma,\partial^2\gamma,\dots]$. Equivalently, $\H^*_{S^1}(\C)_+$ may be described as the vertex algebra ideal generated by $\partial\gamma$. Thus Theorem \circleloc~gives a complete description of $\H^*_{S^1}(\cQ'(M))$ in terms of classical data. 
 \others{Example}{$M = \C\P^1$, where $S^1$ has isotropy weights $-1,1$ at the fixed points $p_0$, $p_1$.}

Since $M^{S^1} = \{p_0,p_1\}$, Theorem \circleloc~shows that \eqn\cpi{\H^*_{S^1}(\cQ'(M))_+ = \H^*_{S^1}(\C)_+\otimes H^*(\{p_0,p_1\}) = \H^*_{S^1}(\C)_+\oplus \H^*_{S^1}(\C)_+.}
It follows that $\H^*_{S^1}(\cQ'(M))$ is the free abelian vertex algebra generated by $H^*_{S^1}(M)$.

\others{Example}{$M = \C\P^2$, where $S^1$ acts with isotropy weights $-1,0,1$ at the fixed points.} We claim that $\H^*_{S^1}(\cQ'(M))$ is {\it not} the vertex algebra $\bra H^*_{S^1}(M)\ket$ generated by the weight zero component.
Classically, $H^*_{S^1}(M)\cong \C[t,\omega]/\bra \omega(\omega-t)(\omega+t)\ket$, where $\omega$ is the equivariant symplectic form and $t$ is the image of the generator of $H^*_{S^1}(pt)$ under the Chern-Weil map. From this description, it is clear that $dim~H^2_{S^1}(M) = 2$. Hence in the vertex subalgebra $\bra H^*_{S^1}(M)\ket \subset \H^*_{S^1}(\cQ'(M))$, the subspace of degree 2 and weight 1 can have dimension at most 2. On the other hand, $M^{S^1}$ consists of three isolated fixed points, so $\H^*_{S^1}(\cQ'(M))_+ = \H^*_{S^1}(\C)_+ \oplus \H^*_{S^1}(\C)_+ \oplus \H^*_{S^1}(\C)_+$ by Theorem \circleloc. In particular,  $dim~\H^2_{S^1}(\cQ'(M))[1] = 3$, so it must contain elements that do not lie in $\bra H^*_{S^1}(M)\ket$.

\theorem{For any $S^1$-manifold $M$, $\H^*_{S^1}(\cQ'(M))$ is generated as a vertex algebra by $\H^*_{S^1}(\cQ'(M))[0]\oplus \H^*_{S^1}(\cQ'(M))[1]$.}\thmlab\generation
\proof Fix a basis $\{\alpha_i|\ i\in I\}$ of $H^*(M^{S^1})$. Consider the collection $$\cC = \{\partial\gamma\otimes \alpha_i|\ i\in I\} \subset \H^*_{S^1}(\C)[1]\otimes H^*(M^{S^1}),$$ and let $\bra \cC\ket$ denote the vertex subalgebra of $\H^*_G(\cQ'(M))$ generated by $\cC$. Since $${\bf i}^*:\H^*_{S^1}(\cQ'(M))\ra\H^*_{S^1}(\cQ'(M^{S^1}))\cong \H^*_{S^1}(\C)\otimes H^*(M^{S^1})$$ preserves $\partial$, it follows that $(\partial^{k+1}\gamma)\otimes\alpha_i=\partial^k (\partial\gamma\otimes\alpha_i) $ lies in $\bra \cC\ket$, for all $k\geq 0$. We claim that $\cC$ together with $H^*_{S^1}(M)$ generates $\H^*_{S^1}(\cQ'(M))$ as a vertex algebra. At weight zero the claim is obvious, so let $\omega\in \H^*_{S^1}(\cQ'(M))_+$ be a nonzero element. In particular, this implies that $M^{S^1}$ is nonempty. 

By Theorem \circleloc, we can write
$\omega = \sum_{i\in I} p_i\otimes \alpha_i$, where $p_i\in \H^*_{S^1}(\C)_+$. Since $p_i$ has positive weight, it is divisible by $\partial^{k_i}\gamma$ for some $k_i>0$, and we may write $p_i = q_i \partial^{k_i}\gamma$ where $q_i\in \H^*_{S^1}(\C)$. Since $M^{S^1}$ is nonempty, the chiral Chern-Weil map is injective, and since $\H^*_{S^1}(\C)$ is generated by $H^*_{S^1}(pt)$ as a vertex algebra, $q_i\otimes 1$ lies in $\bra H^*_{S^1}(M)\ket$. Since $p_i\otimes \alpha_i = ~:(q_i\otimes 1)(\partial^{k_i}\gamma\otimes\alpha_i):$ for each $i\in I$, the claim follows. $\Box$

\corollary{If $M$ is a compact $S^1$-manifold, $\H^*_{S^1}(\cQ'(M))$ is finitely generated as a vertex algebra.}

\subsec{The case $T=S^1\times S^1$ and $M = \C\P^2$}
Next, we consider the case where $T = S^1\times S^1$. This is analogous to the case $G = G_1\times G_2$ with $G_1,G_2$ simple, but can be more subtle because many different copies of $S^1$ inside $T$ can arise as stabilizers of points in $M$. As an example, we compute $\H^*_T(\cQ'(M))$ in the case $M=\C\P^2$ with the usual linear action of $T$. Note that $T$ contains three copies of $S^1$ which arise as stabilizer subgroups, which we denote by $T_i$, $i=1,2,3$. Each $M^{T_i}$ is a copy of $\C\P^1$ which we denote by $M_i$, and $M_i\cap M_j$ consists of a single point $p_{ij}$. Let $U_i$ be a $T$-invariant tubular neighborhood of $M_i$. Clearly $\H^*_T(\cQ'(U_i)) = \H^*_T(\cQ'(M_i)) = \H^*_{T_i}(\C) \otimes \H^*_{T/T_i}(\cQ'(M_i))$, and the action of $T/T_i$ on $M_i$ is the standard action of $S^1$ on $\C\P^1$. 

\lemma{$\H^*_T(\cQ'(U_1\cup U_2))_+$ is linearly isomorphic to 
\eqn\oneandtwo{\H^*_{T_1}(\C)_+\otimes H^*_{T/T_1}(M^{T_1}) \bigoplus \H^*_{T_2}(\C)_+ \otimes H^*_{T/T_2}(M^{T_2})} $$\bigoplus \H^*_{T_1}(\C)_+\otimes \H^*_{T/T_1} (\cQ'(p_{13}))_+\bigoplus \H^*_{T_2}(\C)_+\otimes \H^*_{T/T_2} (\cQ'(p_{23}))_+\bigoplus \H^*_{T_1}(\C)_+\otimes \H^*_{T/T_1}(\cQ'(p_{12}))_+.$$}\thmlab\onetwo

\proof Consider the Mayer-Vietoris sequence
$$\cdots\ra \H^*_T(\cQ'(U_1\cup U_2))_+\ra \H^*_T\cQ'(U_1))_+\oplus \H^*_T(\cQ'(U_2))_+\ra \H^*_T(\cQ'(U_1\cap U_2))_+ \ra \cdots,$$ which we can replace with
\eqn\mayerv{\cdots\ra \H^*_T(\cQ'(U_1\cup U_2))_+\ra \H^*_T(\cQ'(M_1))_+\oplus \H^*_T(\cQ'(M_2))_+\ra \H^*_T(\cQ'(p_{12}))_+ \ra \cdots.}

By \cpi, $\H^*_T(\cQ'(M_1))_+$ is linearly isomorphic to $$\H^*_{T_1}(\C)_+\otimes H^*_{T/T_1}(M_1) \bigoplus H^*_{T_1}(pt)\otimes \H^*_{T/T_1}(\C)_+\otimes H^*(\{p_{12}, p_{13}\}) $$ $$\bigoplus \H^*_{T_1}(\C)_+\otimes \H^*_{T/T_1}(\C)_+\otimes H^*(\{p_{12}, p_{13}\}).$$ Likewise, $\H^*_T(\cQ'(M_2))_+$ is linearly isomorphic to $$\H^*_{T_2}(\C)_+\otimes H^*_{T/T_2}(M_2) \bigoplus H^*_{T_2}(pt)\otimes \H^*_{T/T_2}(\C)_+\otimes H^*(\{p_{12}, p_{23}\})$$ $$ \bigoplus \H^*_{T_2}(\C)_+\otimes \H^*_{T/T_2}(\C)_+\otimes H^*(\{p_{12}, p_{23}\}).$$

Thus the middle term $\H^*_T(\cQ'(M_1))_+\oplus \H^*_T(\cQ'(M_2))_+$ in \mayerv~can be identified with the direct sum of the above six subspaces, so we may write $\omega\in \H^*_T(\cQ'(M_1))_+\oplus \H^*_T(\cQ'(M_2))_+$ as a $6$-tuple $(\omega_1,\dots,\omega_6)$ as in the proof of Theorem \productsimple. The map $$\H^*_T(\cQ'(M_1))_+\oplus \H^*_T(\cQ'(M_2))_+\ra \H^*_T(\cQ'(p_{12}))_+$$ in \mayerv, which we denote by $\phi$, is surjective, so we may identify $\H^*_T(\cQ'(U_1\cup U_2))_+$ with $Ker(\phi)$.

Given $\alpha \in \H^*_{T_1}(\C)_+\otimes H^*_{T/T_1}(M^{T_1})$, we can find $$\tilde{\alpha} \in H^*_{T_2}(pt)\otimes \H^*_{T/T_2}(\C)_+\otimes H^*(\{p_{12}, p_{23}\})$$ such that $(\alpha,0,0,0,\tilde{\alpha},0)\in Ker(\phi)$. The assignment $\alpha\mapsto  (\alpha,0,0,0,\tilde{\alpha},0)$ identifies $ \H^*_{T_1}(\C)_+\otimes H^*_{T/T_1}(M^{T_1})$ with a linear subspace of $Ker(\phi)$. Likewise, given $\alpha \in \H^*_{T_2}(\C)_+\otimes H^*_{T/T_2}(M_2)$, there exists $$\tilde{\alpha}\in H^*_{T_1}(pt)\otimes \H^*_{T/T_1}(\C)_+\otimes H^*(\{p_{12}, p_{13}\})$$ such that $(0,\tilde{\alpha},0,\alpha,0,0)\in Ker(\phi)$, so $ \H^*_{T_2}(\C)_+\otimes H^*_{T/T_2}(M_2)$ may be identified with another subspace of $Ker(\phi)$.

Next, note that $H^*(\{p_{12}, p_{13}\}) = H^*(p_{12})\oplus H^*(p_{13})$, so that $\H^*_{T_1}(\C)_+\otimes \H^*_{T/T_1}(\C)_+\otimes H^*(p_{13})$ may be regarded as a subspace of $\H^*_{T_1}(\C)_+\otimes \H^*_{T/T_1}(\C)_+\otimes H^*(\{p_{12}, p_{13}\})$. Clearly $\H^*_{T_1}(\C)_+\otimes \H^*_{T/T_1}(\C)_+\otimes H^*(p_{13})$ lies in $Ker(\phi)$ since $p_{13}\notin U_1\cap U_2$. Similarly,  $\H^*_{T_2}(\C)_+\otimes \H^*_{T/T_2}(\C)_+\otimes H^*(p_{23})$ is a subspace of  $\H^*_{T_2}(\C)_+\otimes \H^*_{T/T_2}(\C)_+\otimes H^*(\{p_{12}, p_{23}\})$ which lies in $Ker(\phi)$. Finally, note that 
$$\H^*_{T_1}(\C)_+\otimes \H^*_{T/T_1}(\C)_+ =  \H^*_{T_2}(\C)_+\otimes \H^*_{T/T_2}(\C)_+.$$ It follows that $\H^*_{T_1}(\C)_+\otimes \H^*_{T/T_1}(\C)_+\otimes H^*(\{p_{12}, p_{13}\})$ and $\H^*_{T_2}(\C)_+\otimes \H^*_{T/T_2}(\C)_+\otimes H^*(\{p_{12}, p_{23}\})$ each contain a copy of $\H^*_{T_1}(\C)_+\otimes \H^*_{T/T_1}(\C)_+\otimes H^*(p_{12})$. Thus given $\alpha\in \H^*_{T_1}(\C)_+\otimes \H^*_{T/T_1}(\C)_+\otimes H^*(p_{12})$, $(0,0,\alpha,0,0,\alpha)$ will lie in $Ker(\phi)$. This identifies $\H^*_{T_1}(\C)_+\otimes \H^*_{T/T_1}(\C)_+\otimes H^*(p_{12})$ with another subspace of $Ker(\phi)$. Finally, it is straightforward to check that these five subspaces of $Ker(\phi)$ intersect pairwise trivially and account for all of $Ker(\phi)$. This completes the proof of Lemma \onetwo. $\Box$

Next, since $T$ acts locally freely on the complement of $U_1\cup U_2\cup U_3$, 
$$\H^*_T(\cQ'(M))_+ = \H^*_T(\cQ'(U_1\cup U_2\cup U_3))_+,$$ so we have a Mayer-Vietoris sequence
\eqn\mayervii{\cdots\ra \H^*_T(\cQ'(M))_+ \ra \H^*_T(\cQ'(U_1\cup U_2))_+\oplus \H^*_T(\cQ'(U_3))_+\ra\H^*_T(\cQ'((U_1\cup U_2) \cap U_3))_+ .}
Note that $(U_1\cup U_2)\cap U_3=\{p_{13}, p_{23}\}$, so that $$\H^*_T(\cQ'((U_1\cup U_2) \cap U_3))_+  = \H^*_T(\cQ'(p_{13}))_+\oplus \H^*_T(\cQ'(p_{23}))_+ = \H^*_T(\C)_+\oplus \H^*_T(\C)_+.$$ As for the other terms in \mayervii, $\H^*_T(\cQ'(U_1\cup U_2))_+$ is given by \oneandtwo, and $\H^*_T(\cQ'(U_3))_+$ is isomorphic to
$$ \H^*_{T_3}(\C)_+\otimes H^*_{T/T_3}(M_3)\bigoplus H^*_{T_3}(pt) \otimes \H^*_{T/T_3}(\cQ'(M_3))_+\bigoplus \H^*_{T_3}(\C)_+ \otimes \H^*_{T/T_3}(\cQ'(M_3))_+,$$ where 
$$\H^*_{T/T_3}(\cQ'(M_3))_+ = \H^*_{T/T_3}(\cQ'(p_{13}))_+\oplus \H^*_{T/T_3}(\cQ'(p_{23}))_+ = \H^*_{T/T_3}(\C)_+\oplus \H^*_{T/T_3}(\C)_+,$$ by  \cpi. It is easy to check that the map 
$$ \H^*_T(\cQ'(U_1\cup U_2))_+\oplus \H^*_T(\cQ'(U_3))_+\ra\H^*_T(\cQ'((U_1\cup U_2) \cap U_3))_+$$ in \mayervii, which we denote by $\psi$, is surjective, so that $\H^*_T(\cQ'(M))_+ = Ker(\psi)$.

Suppose that $\alpha\in \H^*_{T_3}(\C)_+\otimes H^*_{T/T_3}(M_3)\subset \H^*_T(\cQ'(M_3))_+$. Then there are unique elements $$\alpha_1\in \H^*_{T_1}(\C)_+\otimes H^*_{T/T_1}(M_1)\subset \H^*_T(\cQ'(U_1\cup U_2))_+,$$ $$\alpha_2\in \H^*_{T_2}(\C)_+\otimes H^*_{T/T_2}(M_2)\subset \H^*_T(\cQ'(U_1\cup U_2))_+,$$ such that $(\alpha_1+\alpha_2, \alpha)\in \H^*_T(\cQ'(U_1\cup U_2))_+\oplus \H^*_T(\cQ'(M_3))_+$ lies in $Ker(\psi)$. The assignment $\alpha\mapsto (\alpha_1+\alpha_2, \alpha)$ identifies $\H^*_{T_3}(\C)_+\otimes H^*_{T/T_3}(M_3)$ with a linear subspace of $Ker(\psi)$. Similarly, we may identify each of the spaces 
$$\H^*_{T_1}(\C)_+\otimes H^*_{T/T_1}(M_1)\subset \H^*_T(\cQ'(U_1\cup U_2))_+,$$ $$\H^*_{T_2}(\C)_+\otimes H^*_{T/T_2}(M_2)\subset \H^*_T(\cQ'(U_1\cup U_2))_+,$$
$$\H^*_{T_1}(\C)_+\otimes \H^*_{T/T_1}(\C)_+\otimes H^*(p_{12})\subset \H^*_T(\cQ'(U_1\cup U_2))_+,$$ $$\H^*_{T_2}(\C)_+\otimes \H^*_{T/T_2}(\C)_+\otimes H^*(p_{23})\subset \H^*_T(\cQ'(U_1\cup U_2))_+,$$ $$\H^*_{T_3}(\C)_+\otimes \H^*_{T/T_3}(\C)_+\otimes H^*(p_{13}) \in \H^*_T(\cQ'(M_3))_+,$$ with a subspace of $Ker(\psi)$. It is easy to check that these subspaces intersect pairwise trivially and account for all of $Ker(\psi)$. To summarize, we have proved

\theorem{For $M = \C\P^2$ and $T = S^1\times S^1$ as above, $\H^*_T(\cQ'(M))_+$ is linearly isomorphic to  \eqn\cpii{\big(\bigoplus_{i=1}^3 \H^*_{T_i}(\C)_+ \otimes H^*_{T/T_i}(M_i)\big)\bigoplus \big(\bigoplus_{i=1}^3 \H^*_{T_i}(\C)_+\otimes \H^*_{T/T_i}(\C)_+\big).}}\thmlab\cptwo

As in the case $G = G_1\times G_2$ for $G_1,G_2$ simple, the vertex algebra structure of $\H^*_T(\cQ'(M))$ for $M = \C\P^2$ and $T = S^1\times S^1$ may be deduced from Theorem \cptwo. Let $\cA = \H^*_T(\cQ'(M))$, and $\cB = \H^*_T(\cQ'(U_1\oplus U_2))\oplus \H^*_T(\cQ'(U_3))$, and let $f:\cA\ra \cB$ be the map
$$ \H^*_T(\cQ'(M))\ra  \H^*_T(\cQ'(U_1\oplus U_2))\oplus \H^*_T(\cQ'(U_3))$$ appearing in \mayervii. Clearly the conditions of Lemma \vastruct~are satisfied. As usual, the product $\circ_{-1}$ on $\H^*_T(\cQ'(M))[0] = H^*_T(M)$ is classical, and the vertex algebra structure of $\cB$ is determined by Lemma \onetwo~and the structure of $\H^*_{S^1}(\cQ'(\C\P^1))$, which is given by \cpi. By Lemma \vastruct, the vertex algebra structure of $\H^*_T(\cQ'(M))$ is uniquely determined by this data.

\corollary{For $M = \C\P^2$ and $T = S^1\times S^1$, $\H^*_T(\cQ'(M))$ is generated as a vertex algebra by $\bigoplus_{n = 0}^2 \H^*_T(\cQ'(M))[n].$}

\proof This is immediate from the vertex algebra structure of $\H^*_T(\cQ'(M))$ and the structure of $\H^*_T(\C)$. $\Box$

For a general torus $T$, if $M$ is a $T$-manifold of finite orbit type (i.e., only a finite number of subtori $T'\subset T$ can occur as isotropy groups for points in $M$), we expect that $\H^*_T(\cQ'(M))$ will be generated as a vertex algebra by $\bigoplus_{n=0}^N\H^*_T(\cQ'(M))[n]$ for some $N$. Note that a similar statement for $\H^*_G(\cQ'(M))$ when $G$ is nonabelian is out of reach because it is not known if $\H^*_G(\C)$ is a finitely generated vertex algebra.

\newsec{Concluding Remarks and Open Questions}
Our descriptions of $\H^*_G(\cQ'(M))$ and $\H^*_G(\cQ(M))$ are given relative to the family of vertex algebras $\H^*_K(\C)$ for various connected normal subgroups $K\subset G$. An important open question in this theory is to describe $\H^*_G(\C)$ for any $G$. Note that for $\gg=\gg_1\oplus\cdots\oplus \gg_n$, we have $\H^*_G(\C) = \H^*_{G_1}(\C)\otimes\cdots\otimes\H^*_{G_n}(\C)$, where $G_1,\dots,G_n$ are compact, connected Lie groups with Lie algebras $\gg_1,\dots,\gg_n$, respectively. We already know how to describe $\H^*_G(\C)$ when $G$ is abelian, so it suffices to assume that $G$ is simple.

\others{Question}{Recall from \LLI~that for simple $G$, the weight one subspace $\H^*_G(\C)[1]$ is isomorphic to $Hom_G(\gg,S(\gg^*))$, which is a finitely generated, free module over $\H^*_G(\C)[0] = S(\gg^*)^G$, by a theorem of Kostant. Is there a similar representation-theoretic description of $\H^*_G(\C)[n]$ for any $n$?}

\others{Question}{Is $\H^*_G(\C)$ finitely generated as a vertex algebra? Can we find a set of generators?}

\others{Question}{Can we compute the character $\chi(G) = \sum_{p,n} dim~\H^p_G(\C)[n]~z^p q^n$? Does $\chi(G)$ have any nice properties (modularity, relations to other objects from classical Lie theory, etc.)? For compact $M$, how is $\chi(G,M) = \sum_{p,n} dim~\H^p_G(\cQ'(M))[n]~z^p q^n$ related to other known invariants of $M$?}

\others{Question}{Is there a vertex algebra valued equivariant cohomology $\H^*_G(M)$ for any topological $G$-space $M$ which coincides with $\H^*_G(\cQ'(M))$ when $M$ is a smooth $G$-manifold? In particular, we must have $\H^*_G(pt) = \H^*_G(\C)$; what is the topological interpretation of $\H^*_G(\C)$?}

\footatend\vfill\supereject\immediate\closeout\rfile\writestoppt
\baselineskip=14pt\centerline{{\bf References}}\bigskip{\frenchspacing%
\parindent=20pt\escapechar=` \input refs.tmp\vfill\eject}\nonfrenchspacing

\bs
\item{} Bong H. Lian, Department of Mathematics, Brandeis University, Waltham, MA 02454.

\bs
\item{} Andrew R. Linshaw, Department of Mathematics, University of California, San Diego, La Jolla, CA 92093.

\bs
\item{} Bailin Song, Department of Mathematics, University of California, Los Angeles, Los Angeles, CA 90095.

\end